# CRESTED PRODUCTS OF MARKOV CHAINS

By Daniele D'Angeli and Alfredo Donno

*Université de Genève*

In this work we define two kinds of crested product for reversible Markov chains, which naturally appear as a generalization of the case of crossed and nested product, as in association schemes theory, even if we do a construction that seems to be more general and simple. Although the crossed and nested product are inspired by the study of Gelfand pairs associated with the direct and the wreath product of two groups, the crested products are a more general construction, independent from the Gelfand pairs theory, for which a complete spectral theory is developed. Moreover, the $k$-step transition probability is given. It is remarkable that these Markov chains describe some classical models (Ehrenfest diffusion model, Bernoulli–Laplace diffusion model, exclusion model) and give some generalization of them.

As a particular case of nested product, one gets the classical Insect Markov chain on the ultrametric space. Finally, in the context of the second crested product, we present a generalization of this Markov chain to the case of many insects and give the corresponding spectral decomposition.

**1. Introduction.** The starting point for this work is the paper [4] on Gelfand pairs (for a general theory and applications see [7]) by Ceccherini-Silberstein, Scarabotti and Tolli. In that paper, particular cases of homogeneous spaces are presented. Given two finite groups $G$ and $F$, together with two subgroups $K \leq G$ and $H \leq F$, let $X = G/K$ and $Y = F/H$ be the corresponding homogeneous spaces. The Gelfand pairs associated with the direct product $G \times F$ and the wreath product $F \wr G$ are studied. A complete spectral decomposition, together with the associated spherical functions, is given for the space $L(X \times Y)$.

One can ask if it is possible to get the same spectral decomposition without using the group structure, but considering the product (in some sense) of two reversible Markov chains $P$ and $Q$ defined on two finite sets $X$ and $Y$, respectively.

The analogue of direct product and wreath product is developed in Sections 5 and 6, respectively. In particular, in Section 6, the Markov chain on the ultrametric space, studied by Figà-Talamanca in [11] and that we call "Insect Markov chain," is revisited from a different point of view. See also [6].

An equivalent approach to Gelfand pairs product comes from the association schemes theory (see [2]). In particular, the direct and the wreath product of groups correspond to the crossed and nested product of association schemes. In [3] Bailey and Cameron introduced a more general product for association schemes called the crested product.

Inspired by this definition, in this paper two new products of Markov chains are introduced. We will refer to them as the *first* and the *second* crested product, respectively. Actually, the definition does not seem to be the exact analogue, in terms of Markov chains, of the crested product for association schemes, but it looks even more interesting. In fact, they extend the crossed and the nested product and a complete spectral theory is presented in Sections 4 and 7.

In particular, in Section 7, we show that the second crested product is a generalization of a Markov chain that we will call multi-insect, which is an analogue of the so-called nonbinary Johnson scheme, already studied in [4], that corresponds to considering $h$ insects living in different subtrees and moving only one per each step in such a way that their distance is preserved. We give a generalization of the results of [4] in Section 7 with a more elementary and direct proof. This problem has been suggested by Scarabotti. As an example, the case of the bi-insect is studied.

**2. Basic examples and motivations.** In this section we present some examples, classical and new, of finite diffusion processes that one can treat using methods from harmonic analysis. In the next sections, we will study them by introducing special products of finite Markov chains and performing a complete spectral analysis of the associated Markov operator.

2.1. *Classical models.*

2.1.1. *The Ehrenfest diffusion model.* The classical Ehrenfest diffusion model consists of two urns numbered $0,1$ and $n$ balls numbered $1,\ldots,n$. A configuration is given by a placement of the balls into the urns. Therefore there are $2^n$ configurations and each of them can be identified with the subset $A$ of $\{1,\ldots,n\}$ corresponding to the set of balls contained, for instance, in the urn 0. Note that there is no ordering inside the urns. The initial configuration corresponds to the situation when all balls are inside the urn 0.

Then, at each step, a ball is randomly chosen (with probability $1/n$) and it is moved to the other urn. Denoting by $Q_n$ the set of all subsets of $\{1, \ldots, n\}$, we can describe this random process as follows. If we are in a configuration $A \in Q_n$, then, at the next step, we are in a new configuration $B \in Q_n$ of the form $A \coprod \{j\}$ for some $j \notin A$ or $A \setminus \{j'\}$ for some $j' \in A$ and each of these configurations is chosen with probability $1/n$.

In [5], the authors give a generalization of this model to the case of $m$ urns and $n$ balls.

2.1.2. *The Bernoulli–Laplace diffusion model.* The Bernoulli–Laplace diffusion model consists of two urns numbered $0, 1$ and $2n$ balls, numbered $1, \ldots, 2n$. A configuration is a placement of the balls into the two urns, $n$ balls each. Therefore there are $\binom{2n}{n}$ configurations; each of them can be identified with an $n$-subset $A$ of $\{1, \ldots, 2n\}$ corresponding to the set of balls contained in the urn $0$. The initial configuration corresponds to the situation when the balls contained in the urn $0$ are $1, 2, \ldots, n$. Denoting by $\Omega_n$ the set of all $n$-subsets of $\{1, 2, \ldots, 2n\}$, we can describe this random process as follows. If we are in a configuration $A \in \Omega_n$, then at the next step we are in a new configuration $B \in \Omega_n$ of the form $A \coprod \{i\} \setminus \{j\}$, for some $i \notin A$ and $j \in A$ and each of these configurations is chosen with probability $1/n^2$.

2.1.3. *The exclusion model.* Let $(X, E)$ be an undirected connected graph with $n$ vertices. Fix $h \leq n$ and choose an $h$-subset $A$ of $X$. We put a particle on each vertex of $A$. Then the exclusion process is defined as a Markov chain on the set $\Omega_h$ of $h$-subsets of $X$. Starting from the state $A$, we pick a particle on some vertex of $A$ with probability proportional to the degree of the vertex that it occupies, pick a neighboring site of this vertex at random and move the particle to the neighboring site provided this site is unoccupied. If the site is occupied by another particle, the chain stays at $A$. See [8] and [12] for more examples and details.

2.1.4. *The Insect Markov chain.* This Markov chain, in its classical version introduced by Figà-Talamanca [11], is performed on the (finite) ultrametric space given by the $n$th level of a $q$-ary rooted tree. At the starting instant the insect stays in a fixed vertex of the $n$th level. It performs a simple random walk on the tree. The next state of the Markov chain is given when the insect reaches again the $n$th level. The probability transition between two vertices only depends on their ultrametric distance and so this Markov chain is invariant under the action of the full automorphisms group of the rooted tree.

2.2. *Generalizations.*

2.2.1. *The first crested product.* The Ehrenfest model, generalized to the case of $n$ balls and $m$ urns, is described by the *crossed product* of $n$ Markov chains defined on a space of cardinality $m$ (the space of the urns), as will be shown in Section 5. Actually, the definition of crossed product that we give admits that the choice of the $i$th ball is weighted by a not necessarily homogeneous probability distribution $p_i^0$.

One can ask what happens if a hierarchy in the set of balls is introduced. So suppose to have $n$ ordered balls, numbered by $1, \ldots, n$, and $m$ urns.

At each step, one ball (say the $i$th one) is randomly chosen and it is moved into another urn following a transition probability $P_i$ defined on the space of the urns. For all $j > i$, the ball $j$ is also randomly moved to another urn.

This model is described by the *nested product* of $n$ Markov chains defined on a space of cardinality $m$ (the space of the urns), as will be shown in Section 6, where we suppose that each ball is chosen according with a probability distribution $p_i^0$.

The Insect Markov chain described in Section 6.3 is obtained from this model for a particular value of $p_i^0$ and forcing the transition probability $P_i$ to be uniform, as shown in Proposition 6.4.

The *first crested product* introduced in Section 4 is a generalization of both crossed and nested product and it describes the following diffusion model.

*The $(C, N)$-Ehrenfest model.* We have $n$ balls numbered by $1, \ldots, n$ and $m$ urns. Suppose that we have a partition of the set of the balls in two subsets $C \coprod N$. The balls in $C$ are white, the balls in $N$ are black. At each step, we choose a ball $i$ according with a probability distribution $p_i^0$: if it is white, then we move it to another urn following a transition probability $P_i$ and all the other balls are not moved. If it is black, we move it by $P_i$ and then each ball (white or black) numbered by $j$, for $j > i$, is moved uniformly to a new urn.

*The card players model.* The previous model can also be described in the following way. Suppose that $n$ card players sit at the same edge of a table. Each of them has a deck of cards consisting of $m$ cards. Suppose that the set of players is partitioned in two subsets: blond hair and brown hair. At the starting moment, each player has a card. A player is randomly chosen, so he puts his card inside the deck and then chooses another card. If he is blond, the others keep their own card; if he is brown, then all the players (blond or brown) sitting at his right side also must change their card.

2.2.2. *The second crested product.*

*The bi-insect.* This is a generalization of the Insect Markov chain. In this case we have two insects that, at the starting time, live in the $n$th level of a tree, following the condition that their (ultrametric) distance is maximal, which is equivalent to requiring that they live in two different subtrees of depth $n-1$. At each step, with probability $p_0$ they do not leave their own subtrees; with probability $1-p_0$ one of them can leave its subtree choosing another subtree not occupied by the other insect. In the first case, we randomly choose one insect and change its position according with the Insect Markov chain on the corresponding subtree. In the second case, the moved insect randomly chooses its new position in the new subtree that it has occupied. This yields a Markov chain on the space of all possible configurations of two insects at maximal distance in a tree of depth $n$; this Markov chain will be studied in Section 7.2.

*The multi-insect.* An immediate generalization of the previous construction can be performed considering $h$ insects living in the $n$th level of the tree and having maximal distance from each other. At each step, with probability $p_0$ only one insect moves to another site staying in the same subtree; with probability $1-p_0$ only one insect occupies a new subtree of depth $n-1$ and randomly chooses its position. This situation will be referred to as the *multi-insect* in Section 7.

Moreover, the multi-insect can be also regarded as an exclusion process on the ultrametric space with the constraint that the $h$ insects keep maximal distance from each other. At each step an insect moves, as in the exclusion process a particle does.

The aim of this work is also to diagonalize a more general Markov chain, that we call the *second crested product*, described in Section 7 and that corresponds to the following model.

*The crested Bernoulli–Laplace model.* Suppose that we have two urns numbered $0, 1$ and $n$ balls: the urn 0 contains $h$ balls and the urn 1 contains the remaining $n - h$ balls. The $h$ balls in the urn 0 are endowed with a label that can go from 1 to $s$ (note that two different balls can have the same label). Moreover, we can use a (nonhomogeneous) coin. At each step, we flip the coin: if we get heads, then we randomly choose a ball $i$ in the urn 0 and a ball $j$ in the urn 1 and exchange them. Before putting $j$ in the urn 0 we attach a label (randomly chosen in $\{1, \ldots, s\}$) on it.

If we get tails, then we randomly choose a ball $i$ in the urn 0 and we change or not its label according with a transition probability $Q$ on the set $\{1, \ldots, s\}$.

Observe that, if we get heads, then we perform on the set of the balls an analogue of the Bernoulli–Laplace diffusion model, where the urns contain $h$ and $n - h$ balls, respectively.

Although in Section 7.1 the spectral analysis of the second crested product is given, in Section 7.2 we give explicit computations only in the case of the bi-insect, which seems to be easier to handle.

2.3. *Motivations.* The classical models described in the previous sections are well known in literature and they have been studied by many authors using approaches deriving from Gelfand pairs theory, representations theory, associations schemes, distance regular graphs and orthogonal polynomials.

In the Ehrenfest model, generalized to the case of $m$ urns and $n$ balls, the set of all configurations can be identified with the Hamming scheme $X_{n,m} = \{0, \ldots, m-1\}^n$ consisting of all the functions defined on the set $\{1, \ldots, n\}$ taking values in $\{0, \ldots, m-1\}$. The group that naturally acts on $X_{n,m}$ is the wreath product $S_m \wr S_n$. The stabilizer of a fixed function in $X_{n,m}$ is isomorphic to the subgroup $S_{m-1} \wr S_n$. In [5], for example, the authors show that $(S_m \wr S_n, S_{m-1} \wr S_n)$ is a Gelfand pair, they give the decomposition of the space $L(X_{n,m})$ of the complex functions defined on $X_{n,m}$ into irreducible submodules and provide an expression for the respective spherical functions by using Krawtchouk polynomials. These submodules coincide with the eigenspaces of the Markov operator on $L(X_{n,m})$ describing the Ehrenfest model. A different approach to the Hamming scheme is given in [2] by using association schemes theory.

Following the idea of Diaconis, it is possible to study the rate of convergence to the stationary distribution of finite Markov chains by using Gelfand pairs theory. In particular, given a Markov chain which is invariant under the action of a group $G$, its transition operator can be expressed via a Fourier series where the classical exponentials are replaced by the irreducible representations of $G$. In [10] the authors apply this argument to the generalized Bernoulli–Laplace diffusion model, consisting of two urns containing $h$ and $n-h$ balls, respectively. All the states of this Markov chain constitute the so-called Johnson scheme associated with the Gelfand pair $(S_n, S_h \times S_{n-h})$. The corresponding decomposition into irreducible subrepresentations (and so into eigenspaces for the associated Markov operator) is given, for instance, in [5], where the spherical functions are described in terms of Hahn polynomials. Like the Hamming scheme, also the Johnson scheme can be presented as an association scheme (see [2]).

Finally, the Insect Markov chain can be studied using a Gelfand pair associated with the automorphisms group of the rooted tree (see [6] and [11]). The Fourier analysis that one can develop allows to give a complete spectral analysis and to investigate some ergodic properties.

The approach proposed in this paper is more elementary. We define and study several Markov chains, generalizing the classical models described above, by using only techniques and arguments coming from linear and multilinear algebra. This allows to give directly the spectral analysis of the associated Markov operator, without invoking the theory of groups and their representations. Our construction shows what really are the structures to the base of the analysis of these models and suggests the way to approach more general cases.

In the case of the first crested product we consider the Cartesian product $X = X_1 \times \cdots \times X_n$ of $n$ finite spaces and a Markov operator $P$ on $L(X)$ given by a convex combination of Markov operators on $L(X_i)$ depending on a partition of the set $\{1,\ldots,n\}$. The eigenspaces of the operator $P$ are expressed as certain combinations of the tensor products of the eigenspaces of the Markov operators on every space $L(X_i)$. Their form is easier in the case of the crossed product when none of the spaces $X_i$ is distinguished from the other ones (see Section 5). On the other hand, it becomes more complicated if one introduces a hierarchy. The extremal case gives the nested product (see Section 6), where the space $X$ can be identified with the $n$th level of a rooted tree.

In the second crested product we consider the Cartesian product of two finite spaces $X$ and $Y$, where we are given the eigenspaces of a Markov operator $Q$ defined on $L(Y)$. Set $\Theta_h$ the space of the functions with domain an $h$-subset of $X$ and values in $Y$. We define a Markov operator $P$ [see (6)] on $L(\Theta_h)$ and we investigate its spectral decomposition by using two differential operators $D$ and $D^*$ (see Definition 7.1) and the known spectral analysis of $Q$.

Although the computational aspect can be more or less complicated depending on the singular cases, the methods used do not require anything more than linear or multilinear algebra.

**3. Preliminaries.** The following topics about finite Markov chains can be found in [6].

Consider a finite set $X$, with $|X| = m$. Let $P$ be a stochastic matrix of size $m$ whose rows and columns are indexed by the elements of $X$, so that

$$\sum_{x \in X} p(x_0, x) = 1,$$

for every $x_0 \in X$. Consider the Markov chain on $X$ with transition matrix $P$.

DEFINITION 3.1. The Markov chain $P$ is reversible if there exists a strict probability measure $\pi$ on $X$ such that

$$\pi(x)p(x,y) = \pi(y)p(y,x),$$

for all $x, y \in X$.

We will say that $P$ and $\pi$ are in detailed balance. For a complete treatment about these and related topics see [1].

Define on $L(X) = \{f : X \longrightarrow \mathbb{C}\}$ a scalar product in the following way:
$$\langle f_1, f_2 \rangle_\pi = \sum_{x \in X} f_1(x) \overline{f_2(x)} \pi(x),$$
for all $f_1, f_2 \in L(X)$ and the linear operator $P : L(X) \longrightarrow L(X)$ by
$$(Pf)(x) = \sum_{y \in X} p(x,y) f(y).$$

It is easy to verify that $\pi$ and $P$ are in detailed balance if and only if $P$ is self-adjoint with respect to the scalar product $\langle \cdot, \cdot \rangle_\pi$. Under these hypotheses, it is known that the matrix $P$ can be diagonalized over the reals. Moreover, 1 is always an eigenvalue of $P$ and for any eigenvalue $\lambda$ one has $|\lambda| \leq 1$.

Let $\lambda_z$ be the eigenvalues of the matrix $P$, for every $z \in X$, with $\lambda_{z_0} = 1$. Then there exists an invertible unitary real matrix $U = (u(x,y))_{x,y \in X}$ such that $PU = U\Delta$, where $\Delta = (\lambda_x \delta_x(y))_{x,y \in X}$ is the diagonal matrix whose entries are the eigenvalues of $P$. This equation gives, for all $x, z \in X$,

$$\sum_{y \in X} p(x,y) u(y,z) = u(x,z) \lambda_z. \tag{1}$$

Moreover, we have $U^T D U = I$, where $D = (\pi(x) \delta_x(y))_{x,y \in X}$ is the diagonal matrix of coefficients of $\pi$. This second equation gives, for all $y, z \in X$,

$$\sum_{x \in X} u(x,y) u(x,z) \pi(x) = \delta_y(z). \tag{2}$$

Hence, the first equation tells us that each column of $U$ is an eigenvector of $P$; the second one tells us that these columns are orthogonal with respect to the product $\langle \cdot, \cdot \rangle_\pi$.

PROPOSITION 3.2. *The $k$th step transition probability is given by*
$$p^{(k)}(x,y) = \pi(y) \sum_{z \in X} u(x,z) \lambda_z^k u(y,z), \tag{3}$$
*for all $x, y \in X$.*

PROOF. The proof is a consequence of (1) and (2). In fact, the matrix $U^T D$ is the inverse of $U$, so that $UU^T D = I$. In formulæ, we have
$$\sum_{y \in X} u(x,y) u(z,y) = \frac{1}{\pi(z)} \Delta_z(x).$$
From the equation $PU = U\Delta$ we get $P = U\Delta U^T D$, which gives
$$p(x,y) = \pi(y) \sum_{z \in X} u(x,z) \lambda_z u(y,z).$$

Iterating this argument, we get

$$P^k = U\Delta^k U^T D,$$

which is the assertion. $\square$

Recall that there exists a correspondence between reversible Markov chains and weighted graphs.

DEFINITION 3.3. A weight on a graph $\mathcal{G} = (X, E)$ is a function $w \colon X \times X \longrightarrow [0, +\infty)$ such that:

(1) $w(x, y) = w(y, x)$;
(2) $w(x, y) > 0$ if and only if $x \sim y$.

If $\mathcal{G}$ is a weighted graph, it is possible to associate with $w$ a stochastic matrix $P = (P(x, y))_{x,y \in X}$ on $X$ by setting

$$p(x, y) = \frac{w(x, y)}{W(x)},$$

with $W(x) = \sum_{z \in X} w(x, z)$. The corresponding Markov chain is called the random walk on $\mathcal{G}$. It is easy to prove that the matrix $P$ is in detailed balance with the distribution $\pi$ defined, for every $x \in X$, as

$$\pi(x) = \frac{W(x)}{W},$$

with $W = \sum_{z \in X} W(z)$. Moreover, $\pi$ is strictly positive if $X$ does not contain isolated vertices. The inverse construction can be done. So, if we have a transition matrix $P$ on $X$ which is in detailed balance with the probability $\pi$, then we can define a weight $w$ as $w(x, y) = \pi(x) p(x, y)$. This definition guarantees the symmetry of $w$ and, by setting $E = \{\{x, y\} : w(x, y) > 0\}$, we get a weighted graph.

There are some important relations between the weighted graph associated with a transition matrix $P$ and its spectrum. In fact, it is easy to prove that the multiplicity of the eigenvalue 1 of $P$ equals the number of connected components of $\mathcal{G}$. Moreover, the following propositions hold.

PROPOSITION 3.4. Let $\mathcal{G} = (X, E, w)$ be a finite connected weighted graph and denote $P$ the corresponding transition matrix. Then the following are equivalent:

(1) $\mathcal{G}$ is bipartite;
(2) the spectrum $\sigma(P)$ is symmetric;
(3) $-1 \in \sigma(P)$.

DEFINITION 3.5. Let $P$ be a stochastic matrix. $P$ is ergodic if there exists $n_0 \in \mathbb{N}$ such that
$$p^{(n_0)}(x,y) > 0 \qquad \text{for all } x, y \in X.$$

PROPOSITION 3.6. Let $\mathcal{G} = (X, E)$ be a finite graph. Then the following conditions are equivalent:

(1) $\mathcal{G}$ is connected and not bipartite;
(2) for every weight function on $X$, the associated transition matrix $P$ is ergodic.

So we can conclude that a reversible transition matrix $P$ is ergodic if and only if the eigenvalue 1 has multiplicity 1 and $-1$ is not an eigenvalue.

In what follows we always suppose that the eigenvalue 1 has multiplicity 1, so that the graph associated with the probability $P$ is connected. This is equivalent to requiring that the probability $P$ is irreducible, according with the following definition.

DEFINITION 3.7. A stochastic matrix $P$ on a set $X$ is irreducible if, for every $x_1, x_2 \in X$, there exists $n = n(x_1, x_2)$ such that $p^{(n)}(x_1, x_2) > 0$.

**4. The first crested product.** In this section we introduce a particular product of Markov chains defined on different sets. This idea is inspired by the definition of crested product for association schemes given in [3].

Let $X_i$ be a finite set, with $|X_i| = m_i$, for every $i = 1, \ldots, n$, so that we can identify $X_i$ with the set $\{0, 1, \ldots, m_i - 1\}$. Let $P_i$ be an irreducible Markov chain on $X_i$ and let $p_i$ be the transition probability associated with $P_i$. Moreover, assume that $p_i$ is in detailed balance with the strict probability measure $\sigma_i$ on $X_i$, that is,
$$\sigma_i(x) p_i(x, y) = \sigma_i(y) p_i(y, x),$$
for all $x, y \in X_i$.

Consider the product $X_1 \times \cdots \times X_n$. Let $\{1, \ldots, n\} = C \coprod N$ be a partition of the set $\{1, \ldots, n\}$ and let $p_1^0, p_2^0, \ldots, p_n^0$ be a probability distribution on $\{1, \ldots, n\}$, that is, $p_i^0 > 0$ for every $i = 1, \ldots, n$ and $\sum_{i=1}^n p_i^0 = 1$.

DEFINITION 4.1. The first crested product of Markov chains $P_i$'s with respect to the partition $\{1, \ldots, n\} = C \coprod N$ is the Markov chain on the product $X_1 \times \cdots \times X_n$ whose transition matrix is
$$P = \sum_{i \in C} p_i^0 (I_1 \otimes \cdots \otimes I_{i-1} \otimes P_i \otimes I_{i+1} \otimes \cdots \otimes I_n)$$
$$+ \sum_{i \in N} p_i^0 (I_1 \otimes \cdots \otimes I_{i-1} \otimes P_i \otimes J_{i+1} \otimes \cdots \otimes J_n),$$

where $I_i$ denotes the identity matrix of size $m_i$ and $J_i$ denotes the uniform matrix on $X_i$, that is,

$$J_i = \frac{1}{m_i}\begin{pmatrix} 1 & 1 & \cdots & 1 \\ 1 & \ddots & & \vdots \\ \vdots & & \ddots & \vdots \\ 1 & \cdots & \cdots & 1 \end{pmatrix}.$$

In other words, we choose an index $i$ in $\{1,\ldots,n\}$ following the distribution $p_1^0,\ldots,p_n^0$. If $i \in C$, then $P$ acts on the $i$th coordinate by the matrix $P_i$ and fixes the remaining coordinates; if $i \in N$, then $P$ fixes the coordinates $\{1,\ldots,i-1\}$, acts on the $i$th coordinate by the matrix $P_i$ and changes uniformly the remaining ones.

For all $(x_1,\ldots,x_n),(y_1,\ldots,y_n) \in X_1 \times \cdots \times X_n$, the transition probability $p$ associated with $P$ is given by

$$p((x_1,\ldots,x_n),(y_1,\ldots,y_n))$$
$$= \sum_{i \in C} p_i^0 (\delta_1(x_1,y_1) \cdots \delta_{i-1}(x_{i-1},y_{i-1})$$
$$\times p_i(x_i,y_i)\delta_{i+1}(x_{i+1},y_{i+1}) \cdots \delta_n(x_n,y_n))$$
$$+ \sum_{i \in N} p_i^0 \left( \frac{\delta_1(x_1,y_1) \cdots \delta_{i-1}(x_{i-1},y_{i-1}) p_i(x_i,y_i)}{\prod_{j=i+1}^n m_j} \right),$$

where $\delta_i$ is defined by

$$\delta_i(x_i,y_i) = \begin{cases} 1, & \text{if } x_i = y_i, \\ 0, & \text{otherwise.} \end{cases}$$

We want to study the spectral theory of the operator $P$. We recall that the following isomorphism holds:

$$L(X_1 \times \cdots \times X_n) \cong \bigotimes_{i=1}^n L(X_i),$$

with $(f_1 \otimes \cdots \otimes f_n)(x_1,\ldots,x_n) := f_1(x_1)f_2(x_2) \cdots f_n(x_n)$.

Assume that, for every $i = 1,\ldots,n$, the following spectral decomposition holds:

$$L(X_i) = \bigoplus_{j_i=0}^{r_i} V_{j_i}^i,$$

that is, $V_{j_i}$ is an eigenspace for $P_i$ with associated eigenvalue $\lambda_{j_i}$ and whose dimension is $m_{j_i}$.

Now set $N = \{i_1,\ldots,i_l\}$ and $C = \{c_1,\ldots,c_h\}$, with $h + l = n$ and such that $i_1 < \cdots < i_l$ and $c_1 < \cdots < c_h$.

THEOREM 4.2. *The probability $P$ defined above is reversible if and only if $P_k$ is symmetric for every $k > i_1$. If this is the case, $P$ is in detailed balance with the strict probability measure $\pi$ on $X_1 \times \cdots \times X_n$ given by*

$$\pi(x_1,\ldots,x_n) = \frac{\sigma_1(x_1)\sigma_2(x_2)\cdots\sigma_{i_1}(x_{i_1})}{m_{i_1+1}\cdots m_n}.$$

PROOF. Consider the elements $x = (x_1,\ldots,x_n)$ and $y = (y_1,\ldots,y_n)$ belonging to $X_1 \times \cdots \times X_n$. First, we want to prove that the condition $\sigma_k = \frac{1}{m_k}$, for every $k > i_1$, is sufficient. Let $k \in \{1,\ldots,n\}$ such that $x_i = y_i$ for every $i = 1,\ldots,k-1$ and $x_k \neq y_k$. Suppose $k < i_1$. Then we have

$$p(x,y) = p_k^0(p_k(x_k,y_k)\delta_{k+1}(x_{k+1},y_{k+1})\cdots\delta_n(x_n,y_n)).$$

If $x_i = y_i$ for every $i = k+1,\ldots,n$, we get

$$\pi(x)p(x,y) = \sigma_1(x_1)\cdots\sigma_k(x_k)\cdots\sigma_{i_1}(x_{i_1})p_k^0\frac{p_k(x_k,y_k)}{m_{i_1+1}\cdots m_n}$$

$$= \sigma_1(y_1)\cdots\sigma_k(y_k)\cdots\sigma_{i_1}(y_{i_1})p_k^0\frac{p_k(y_k,x_k)}{m_{i_1+1}\cdots m_n}$$

$$= \pi(y)p(y,x),$$

since $\sigma_k(x_k)p_k(x_k,y_k) = \sigma_k(y_k)p_k(y_k,x_k)$. If the condition $x_i = y_i$ is not satisfied for every $i = k+1,\ldots,n$, then the equality $\pi(x)p(x,y) = \pi(y)p(y,x) = 0$ easily follows.

If $k = i_1$, then we get

$$p(x,y) = p_{i_1}^0\left(p_{i_1}(x_{i_1},y_{i_1})\frac{1}{m_{i_1+1}\cdots m_n}\right)$$

and so

$$\pi(x)p(x,y) = \sigma_1(x_1)\cdots\sigma_{i_1}(x_{i_1})p_{i_1}^0\frac{p_{i_1}(x_{i_1},y_{i_1})}{m_{i_1+1}^2\cdots m_n^2}$$

$$= \sigma_1(y_1)\cdots\cdots\sigma_{i_1}(y_{i_1})p_{i_1}^0\frac{p_{i_1}(y_{i_1},x_{i_1})}{m_{i_1+1}^2\cdots m_n^2}$$

$$= \pi(y)p(y,x),$$

since $\sigma_{i_1}(x_{i_1})p_{i_1}(x_{i_1},y_{i_1}) = \sigma_{i_1}(y_{i_1})p_{i_1}(y_{i_1},x_{i_1})$.

In the case $k > i_1$, we have

$$p(x,y) = \sum_{i\in N, i\leq k} p_i^0\frac{p_i(x_i,y_i)}{m_{i+1}\cdots m_n}$$

and so

$$\pi(x)p(x,y) = \frac{\sigma_1(x_1)\cdots\sigma_{i_1}(x_{i_1})}{m_{i_1+1}\cdots m_n} \sum_{i\in N, i\leq k} p_i^0 \frac{p_i(x_i,y_i)}{m_{i+1}\cdots m_n}$$

$$= \frac{\sigma_1(y_1)\cdots\sigma_{i_1}(y_{i_1})}{m_{i_1+1}\cdots m_n} \sum_{i\in N, i\leq k} p_i^0 \frac{p_i(y_i,x_i)}{m_{i+1}\cdots m_n}$$

$$= \pi(y)p(y,x).$$

In fact, the terms corresponding to an index $i < k$ satisfy $p_i(x_i, y_i) = p_i(y_i, x_i)$ since $x_i = y_i$; the term corresponding to the index $k$ satisfies $p_k(x_k, y_k) = p_k(y_k, x_k)$ since the equality

$$p_k(x_k, y_k) = p_k(y_k, x_k)$$

holds by hypothesis.

Now we want to prove that the condition $\sigma_k = \frac{1}{m_k}$, for every $k > i_1$, is also necessary. Suppose that the equality $\pi(x)p(x,y) = \pi(y)p(y,x)$ holds. By the hypothesis of irreducibility we can consider two elements $x^0, y^0 \in X_1 \times \cdots \times X_n$ such that $x_{i_1}^0 \neq y_{i_1}^0$ and with the property that $p_{i_1}(x_{i_1}^0, y_{i_1}^0) \neq 0$. Now we have

$$\pi(x^0)p(x^0, y^0) = \pi(y^0)p(y^0, x^0) \Leftrightarrow$$
$$\pi(x^0)p_{i_1}(x_{i_1}^0, y_{i_1}^0) = \pi(y^0)p_{i_1}(y_{i_1}^0, x_{i_1}^0).$$

This gives

$$\frac{\pi(x^0)}{\pi(y^0)} = \frac{p_{i_1}(y_{i_1}^0, x_{i_1}^0)}{p_{i_1}(x_{i_1}^0, y_{i_1}^0)} = \frac{\sigma_{i_1}(x_{i_1}^0)}{\sigma_{i_1}(y_{i_1}^0)}.$$

Consider now the element $x = (x_1^0, \ldots, x_{i_1}^0, y_{i_1+1}^0, \ldots, y_n^0)$. The equality $\pi(x)p(x, y^0) = \pi(y^0)p(y^0, x)$ implies

$$\frac{\pi(x)}{\pi(y^0)} = \frac{p_{i_1}(y_{i_1}^0, x_{i_1}^0)}{p_{i_1}(x_{i_1}^0, y_{i_1}^0)} = \frac{\sigma_{i_1}(x_{i_1}^0)}{\sigma_{i_1}(y_{i_1}^0)}.$$

So we get $\pi(x^0) = \pi(x)$, that is, the probability $\pi$ does not depend on the coordinates $i_1 + 1, \ldots, n$. Set now $x' = (x_1^0, \ldots, x_{i_1}^0, \ldots, x_{k-1}^0, x_k, \ldots, x_n)$. The equality $\pi(x^0)p(x^0, x') = \pi(x')p(x', x^0)$ gives

$$\pi(x^0)\left(\sum_{j\in N, j\leq k} p_j^0(p_j(x_j^0, x_j'))\right) = \pi(x')\left(\sum_{j\in N, j\leq k} p_j^0(p_j(x_j', x_j^0))\right).$$

Since the probability $\pi$ does not depend on the coordinates $i_1 + 1, \ldots, n$, we get $p_k(x_k^0, x_k') = p_k(x_k', x_k^0)$. This implies $\sigma_k(x_k') = \sigma_k(x_k^0)$ and so the hypothesis of irreducibility guarantees that $\sigma_k$ is uniform on $X_k$. This completes the proof. $\square$

THEOREM 4.3. *The eigenspaces of the operator $P$ are given by:*

- $W^1 \otimes \cdots \otimes W^{k-1} \otimes V_{j_k}^k \otimes V_0^{k+1} \otimes V_0^{k+2} \otimes \cdots \otimes V_0^n$, with $j_k \neq 0$, for $k \in \{i_1+1,\ldots,n\}$ and where

$$W^i = \begin{cases} L(X_i), & \text{if } i \in N, \\ V_{j_i}^i, \quad j_i = 0,\ldots,r_i, & \text{if } i \in C, \end{cases}$$

with eigenvalue

$$\sum_{i \in C: i < k} p_i^0 \lambda_{j_i} + p_k^0 \lambda_{j_k} + \sum_{i > k} p_i^0.$$

- $V_{j_1}^1 \otimes \cdots \otimes V_{j_{i_1-1}}^{i_1-1} \otimes V_{j_{i_1}}^{i_1} \otimes V_0^{i_1+1} \otimes \cdots \otimes V_0^n$, with $j_t = 0,\ldots,r_t$, for every $t = 1,\ldots,i_1$, with eigenvalue

$$\sum_{i=1}^{i_1} p_i^0 \lambda_{j_i} + \sum_{i=i_1+1}^n p_i^0.$$

PROOF. Fix an index $k \in \{i_1+1, i_1+2,\ldots,n\}$ and consider the function $\varphi$ in the space

$$W^1 \otimes \cdots \otimes W^{k-1} \otimes V_{j_k}^k \otimes V_0^{k+1} \otimes V_0^{k+2} \otimes \cdots \otimes V_0^n,$$

with $j_k \neq 0$ and

$$W^i = \begin{cases} L(X_i), & \text{if } i \in N, \\ V_{j_i}^i, \quad j_i = 0,\ldots,r_i, & \text{if } i \in C, \end{cases}$$

so that $\varphi = \varphi_1 \otimes \cdots \otimes \varphi_{k-1} \otimes \varphi_k \otimes \varphi_{k+1} \otimes \cdots \otimes \varphi_n$ with $\varphi_i \in W^i$ for $i = 1,\ldots,k-1$, $\varphi_k \in V_{j_k}^k$ and $\varphi_l \in V_0^l$ for $l = k+1,\ldots,n$. Set $x = (x_1,\ldots,x_n)$ and $y = (y_1,\ldots,y_n)$; then

$$(P\varphi)(x) = \sum_y p(x,y)\varphi(y)$$

$$= \sum_y \Bigg(\sum_{i \in C} p_i^0 \delta_1(x_1,y_1)\cdots\delta_{i-1}(x_{i-1},y_{i-1})$$

$$\times p_i(x_i,y_i)\delta_{i+1}(x_{i+1},y_{i+1})\cdots\delta_n(x_n,y_n)$$

$$+ \sum_{i \in N} p_i^0 \delta_1(x_1,y_1)\cdots\delta_{i-1}(x_{i-1},y_{i-1})p_i(x_i,y_i)\frac{1}{m_{i+1}}\cdots\frac{1}{m_n}\Bigg)$$

$$\times \varphi_1(y_1)\cdots\varphi_{k-1}(y_{k-1})\varphi_k(y_k)\varphi_{k+1}(y_{k+1})\cdots\varphi_n(y_n)$$

$$= \sum_{i \in C, i \leq k} \Bigg(\sum_{y_i} p_i^0 p_i(x_i,y_i)\varphi_i(y_i)\Bigg)$$

$$\times \varphi_1(x_1)\cdots\varphi_{i-1}(x_{i-1})\varphi_{i+1}(x_{i+1})\cdots\varphi_n(x_n)$$
$$+ \sum_{i\in C, i>k}\left(\sum_{y_i} p_i^0 p_i(x_i, y_i)\varphi_i(y_i)\right)$$
$$\times \varphi_1(x_1)\cdots\varphi_{i-1}(x_{i-1})\varphi_{i+1}(x_{i+1})\cdots\varphi_n(x_n)$$
$$+ \sum_{i\in N, i>k}\left(\sum_{y_i,\ldots,y_n} p_i^0 p_i(x_i, y_i)\frac{1}{m_{i+1}}\cdots\frac{1}{m_n}\varphi_i(y_i)\cdots\varphi_n(y_n)\right)$$
$$\times \varphi_1(x_1)\cdots\varphi_{i-1}(x_{i-1})$$
$$+ \chi_N(k) \sum_{y_k,\ldots,y_n} p_k^0 p_k(x_k, y_k)\frac{1}{m_{k+1}}\cdots\frac{1}{m_n}$$
$$\times \varphi_1(x_1)\cdots\varphi_{k-1}(x_{k-1})\varphi_k(y_k)\cdots\varphi_n(y_n)$$
$$= \sum_{i\in C, i\leq k} p_i^0 \lambda_{j_i}\varphi(x) + \sum_{i\in C, i>k} p_i^0 \cdot 1 \cdot \varphi(x)$$
$$+ \sum_{i\in N, i>k}\left(\sum_{y_i} p_i^0 p_i(x_i, y_i)\varphi_i(y_i)\right)$$
$$\times \varphi_1(x_1)\cdots\varphi_{i-1}(x_{i-1})\varphi_{i+1}(x_{i+1})\cdots\varphi_n(x_n)$$
$$+ \chi_N(k)\sum_{y_k} p_k^0 p_k(x_k, y_k)$$
$$\times \varphi_1(x_1)\cdots\varphi_{k-1}(x_{k-1})\varphi_k(y_k)\varphi_{k+1}(x_{k+1})\cdots\varphi_n(x_n)$$
$$= \sum_{i\in C, i\leq k} p_i^0 \lambda_{j_i}\varphi(x) + \sum_{i\in C, i>k} p_i^0 \varphi(x)$$
$$+ \sum_{i\in N, i>k} p_i^0\varphi(x) + \chi_N(k)p_k^0\lambda_{j_k}\varphi(x)$$
$$= \left(\sum_{i\in C, i<k} p_i^0\lambda_{j_i} + p_k^0\lambda_{j_k} + \sum_{i>k} p_i^0\right)\varphi(x),$$

where $\chi_N$ is the characteristic function of $N$. Note that in this case the addends corresponding to the indices $i < k$, $i \in N$, are equal to 0 since we have supposed $j_k \neq 0$.

Consider now the function $\varphi$ in the space

$$V_{j_1}^1 \otimes \cdots V_{j_{i_1-1}}^{i_1-1} \otimes V_{j_{i_1}}^{i_1} \otimes V_0^{i_1+1} \otimes \cdots \otimes V_0^n,$$

with $j_t = 0,\ldots,r_t$, for every $t = 1,\ldots,i_1$. In this case we have

$$(P\varphi)(x) = \sum_y p(x,y)\varphi(y)$$

$$
\begin{aligned}
&= \sum_{i\in C, i<i_1} \left(\sum_{y_i} p_i^0 p_i(x_i, y_i)\varphi_i(y_i)\right) \\
&\qquad \times \varphi_1(x_1)\cdots\varphi_{i-1}(x_{i-1})\varphi_{i+1}(x_{i+1})\cdots\varphi_n(x_n) \\
&+ \sum_{i\in C, i>i_1} \left(\sum_{y_i} p_i^0 p_i(x_i, y_i)\varphi_i(y_i)\right) \\
&\qquad \times \varphi_1(x_1)\cdots\varphi_{i-1}(x_{i-1})\varphi_{i+1}(x_{i+1})\cdots\varphi_n(x_n) \\
&+ \sum_{i\in N, i>i_1} \left(\sum_{y_i,\ldots,y_n} p_i^0 p_i(x_i, y_i)\frac{1}{m_{i+1}}\cdots\frac{1}{m_n}\varphi_i(y_i)\cdots\varphi_n(x_n)\right) \\
&\qquad \times \varphi_1(x_1)\cdots\varphi_{i-1}(x_{i-1}) \\
&+ \sum_{y_{i_1},\ldots,y_n} \left(p_{i_1}^0 p_{i_1}(x_{i_1}, y_{i_1})\frac{1}{m_{i_1+1}}\cdots\frac{1}{m_n}\varphi_{i_1}(y_{i_1})\cdots\varphi_n(x_n)\right) \\
&\qquad \times \varphi_1(x_1)\cdots\varphi_{i_1-1}(x_{i_1-1}) \\
&= \sum_{i\in C, i<i_1} p_i^0 \lambda_{j_i}\varphi(x) + \sum_{i\in C, i>i_1} p_i^0 \varphi(x) \\
&\quad + \sum_{i\in N, i>i_1} p_i^0 \varphi(x) + p_{i_1}^0 \lambda_{j_{i_1}}\varphi(x) \\
&= \left(\sum_{i=1}^{i_1} p_i^0 \lambda_{j_i} + \sum_{i=i_1+1}^{n} p_i^0\right)\varphi(x).
\end{aligned}
$$

Observe that, by computing the sum of the dimensions of these eigenspaces, we get

$$\sum_{k=i_1+1}^{n} m_1\cdots m_{k-1}(m_k - 1) + m_1 m_2 \cdots m_{i_1} = m_1 m_2 \cdots m_n,$$

which is just the dimension of the space $X_1 \times \cdots \times X_n$. □

REMARK 4.4. The expression of the eigenvalues of $P$ given in the previous theorem tells us that if $P_i$ is ergodic for every $i = 1, \ldots, n$, then also $P$ is ergodic, since the eigenvalue 1 is obtained with multiplicity 1 and the eigenvalue $-1$ can never be obtained.

We can now give the matrices $U, D$ and $\Delta$ associated with $P$. For every $i$, let $U_i$, $D_i$ and $\Delta_i$ be the matrices of eigenvectors, of the coefficients of $\sigma_i$ and of eigenvalues for the probability $P_i$, respectively. The expression of the matrix $U$, whose columns are an orthonormal basis of eigenvectors for $P$, easily follows from Theorem 4.3. In order to get the diagonal matrix $D$, whose entries are the coefficients of $\pi$, it suffices to consider the tensor product of the corresponding matrices associated with the probability $P_i$, for every $i = 1, \ldots, n$, as it follows from Theorem 4.2. Finally, to get the matrix $\Delta$ of eigenvalues of $P$ it suffices to replace, in the expression of the matrix $P$, the matrix $P_i$ by $\Delta_i$ and the matrix $J_i$ by the corresponding diagonal matrix $J_i^{\text{diag}}$, which has the eigenvalue 1 with multiplicity 1 and the eigenvalue 0 with multiplicity $m_i - 1$. So we have the following proposition.

PROPOSITION 4.5. *Let $P$ be the crested product of the Markov chains $P_i$, with $i = 1, \ldots, n$. Then we have:*

- $U = \sum_{k=i_1+1}^{n} M_1 \otimes \cdots \otimes M_{k-1} \otimes (U_k - A_k) \otimes A_{k+1} \otimes \cdots \otimes A_n + U_1 \otimes U_2 \otimes \cdots \otimes U_{i_1} \otimes A_{i_1+1} \otimes \cdots \otimes A_n$, *with*

$$M_i = \begin{cases} I_i^{\sigma_i\text{-norm}}, & \text{if } i \in N, \\ U_i, & \text{if } i \in C, \end{cases}$$

*where*

$$I_i^{\sigma_i\text{-norm}} = \begin{pmatrix} \frac{1}{\sqrt{\sigma_i(0)}} & & & \\ & \frac{1}{\sqrt{\sigma_i(1)}} & & \\ & & \ddots & \\ & & & \frac{1}{\sqrt{\sigma_i(m_i-1)}} \end{pmatrix}.$$

*By $A_i$ we denote the matrix of size $m_i$ whose entries on the first column are all 1 and the remaining ones are 0:*

- $D = \bigotimes_{i=1}^{n} D_i$.
- $\Delta = \sum_{i \in C} p_i^0 (I_1 \otimes \cdots \otimes I_{i-1} \otimes \Delta_i \otimes I_{i+1} \otimes \cdots \otimes I_n) + \sum_{i \in N} p_i^0 (I_1 \otimes \cdots \otimes I_{i-1} \otimes \Delta_i \otimes J_{i+1}^{\text{diag}} \otimes \cdots \otimes J_n^{\text{diag}})$.

Observe that another matrix $U'$ of eigenvectors for $P$ is given by $U' = \bigotimes_{i=1}^{n} U_i$. The matrix $U$ that we have given above seems to be more useful whenever one wants to compute the $k$th step transition probability $p^{(k)}(0, x)$ using the formula (3), since it contains a greater number of zeros in the first row with respect to $U'$ and so a small number of terms in the sum are nonzero.

Suppose $x = (x_1, \ldots, x_n)$ and $y = (y_1, \ldots, y_n)$ elements in $X = X_1 \times \cdots \times X_n$. From (3) and Proposition 4.5, we have

$p^{(k)}(x, y)$

$$= \pi(y) \Bigg[ \sum_{z \in X} \Bigg( \sum_{r=i_1+1}^{n} m_1(x_1, z_1) \cdots m_{r-1}(x_{r-1}, z_{r-1})(u_r - a_r)(x_r, z_r)$$

$$\times a_{r+1}(x_{r+1},z_{r+1})\cdots a_n(x_n,z_n)$$
$$+ u_1(x_1,z_1)\cdots u_{i_1}(x_{i_1},z_{i_1})$$
$$\times a_{i_1+1}(x_{i_1+1},z_{i_1+1})\cdots a_n(x_n,z_n)\bigg)\lambda_z^k$$
$$\times \bigg(\sum_{r=i_1+1}^n m_1(y_1,z_1)\cdots m_{r-1}(y_{r-1},z_{r-1})(u_r - a_r)(y_r,z_r)$$
$$\times a_{r+1}(y_{r+1},z_{r+1})\cdots a_n(y_n,z_n)$$
$$+ u_1(y_1,z_1)\cdots u_{i_1}(y_{i_1},z_{i_1})$$
$$\times a_{i_1+1}(y_{i_1+1},z_{i_1+1})\cdots a_n(y_n,z_n)\bigg)\bigg],$$

where $m_i, u_i, a_i$ are the probabilities associated with the matrices $M_i, U_i, A_i$ occurring in Proposition 4.5.

**5. The crossed product.** The crossed product of the Markov chains $P_i$'s can be obtained as a particular case of the crested product, by setting $C = \{1,\ldots,n\}$ (so that $N = \varnothing$) and it is also called the *direct* product. The analogous case for product of groups has been studied in [9].

In this case, we get the following transition probability:
$$p((x_1,\ldots,x_n),(y_1,\ldots,y_n)) = \sum_{i=1}^n p_i^0 \delta(x_1,y_1)\cdots p_i(x_i,y_i)\cdots \delta(x_n,y_n).$$

This corresponds to choosing the $i$th coordinate with probability $p_i^0$ and to changing it according with the probability transition $P_i$. So we get
$$p((x_1,\ldots,x_n),(y_1,\ldots,y_n)) = \begin{cases} p_i^0 p_i(x_i,y_i), & \text{if } x_j = y_j \text{ for all } j \neq i, \\ 0, & \text{otherwise.} \end{cases}$$

So, for $X_1 = \cdots = X_n =: X$ and $p_0^1 = \cdots = p_n^0 = \frac{1}{n}$, the probability $p$ defines an analogue of the Ehrenfest model, where $n$ is the number of balls and $|X| = m$ is the number of urns. In order to obtain a new configuration, we choose a ball with probability $1/n$ (let it be the $i$th ball in the urn $x_i$) and with probability $p_i(x_i,y_i)$ we put it in the urn $y_i$.

As a consequence of Theorem 4.2, we get that if $P_i$ is in detailed balance with $\pi_i$, then $P$ is in detailed balance with the strict probability measure $\pi$ on $X_1 \times \cdots \times X_n$ defined as
$$\pi(x_1,\ldots,x_n) = \pi_1(x_1)\pi_2(x_2)\cdots \pi_n(x_n).$$

The matrix $P$ associated with the probability $p$ is given by

(4) $$P = \sum_{i=1}^n p_i^0 (I_1 \otimes \cdots \otimes P_i \otimes \cdots \otimes I_n).$$

The following proposition studies the spectral theory of the operator $P$ and it is a straightforward consequence of Theorem 4.3.

PROPOSITION 5.1. *Let $\varphi_0^i = 1_{X_i}, \varphi_1^i, \ldots, \varphi_{m_i-1}^i$ be the eigenfunctions of $P_i$ associated with the eigenvalues $\lambda_0^i = 1, \lambda_1^i, \ldots, \lambda_{m_i-1}^i$, respectively. Then the eigenvalues of $P$ are the $m_1 m_2 \cdots m_n$ numbers*

$$\lambda_I = \sum_{k=1}^n p_k^0 \lambda_{i_k}^k,$$

*with $I = (i_1, \ldots, i_n) \in \{0, \ldots, m_1 - 1\} \times \cdots \times \{0, \ldots, m_n - 1\}$. The corresponding eigenfunctions are defined as*

$$\varphi_I((x_1, \ldots, x_n)) = \varphi_{i_1}^1(x_1) \cdots \varphi_{i_n}^n(x_n).$$

As a consequence of Proposition 4.5, in order to get the matrices $U, D$ and $\Delta$ associated with $P$, it suffices to consider the tensor product of the corresponding matrices associated with the probability $P_i$, for every $i = 1, \ldots, n$. For every $i$, let $U_i$, $D_i$ and $\Delta_i$ be the matrices of eigenvectors, of the coefficients of $\pi_i$ and of eigenvalues for the probability $P_i$, respectively. We have the following corollary.

COROLLARY 5.2. *Let $P$ be the probability defined in (4). Then we have*

$$PU = U\Delta,$$
$$U^T DU = I,$$

*where $U = \bigotimes_{i=1}^n U_i$, $\Delta = \bigotimes_{i=1}^n \Delta_i$ and $D = \bigotimes_{i=1}^n D_i$.*

In particular, we can express the $k$th step transition probability matrix as

$$P^k = \left(\bigotimes_{i=1}^n U_i\right)\left(\bigotimes_{i=1}^n \Delta_i\right)^k \left(\bigotimes_{i=1}^n U_i\right)^T \left(\bigotimes_{i=1}^n D_i\right).$$

Let $x = (x_1, \ldots, x_n)$ and $y = (y_1, \ldots, y_n)$. Then we get

$$p^{(k)}(x, y) = \pi(y) \sum_I \varphi_I(x) \lambda_I^k \varphi_I(y)$$

$$= \pi_1(y_1) \cdots \pi_n(y_n)$$
$$\times \sum_I \varphi_{i_1}^1(x_1) \cdots \varphi_{i_n}^n(x_n)(p_1^0 \lambda_{i_1}^1 + \cdots + p_n^0 \lambda_{i_n}^n)^k \varphi_{i_1}^1(y_1) \cdots \varphi_{i_n}^n(y_n),$$

with $I = (i_1, \ldots, i_n)$.

As we said in Remark 4.4, if the matrix $P_i$ is ergodic for every $i = 1, \ldots, n$, then also the matrix $P$ is ergodic, since the eigenvalue 1 can be obtained only by choosing $I = 0^n$ and the eigenvalue $-1$ can never be obtained.

REMARK 5.3. Put $n = 1$ and set $X = \{0, 1, \ldots, m-1\}$. Consider the action of the symmetric group $S_m$ on $X$. The stabilizer of a fixed element $x_0 = 0$ is isomorphic to the symmetric group $S_{m-1}$. It is well known (see [5]) that $(S_m, S_{m-1})$ is a Gelfand pair and the following decomposition of $L(X)$ into irreducible representations holds:

$$L(X) = V_0 \oplus V_1,$$

where $V_0 \cong \mathbb{C}$ is the space of constant functions on $X$ and $V_1 = \{f : X \longrightarrow \mathbb{C} : \sum_{i=0}^{m-1} f(i) = 0\}$. So we have $\dim V_0 = 1$ and $\dim V_1 = m - 1$.

Analogously, one can consider the action of the wreath product $S_m \wr S_n$ on $X^n = X \times \cdots \times X$, defined in the natural way, and then one can study the decomposition of $L(X^n)$. We have

$$L(X^n) \cong L(X)^{\otimes n} \cong \bigoplus_{j=0}^{n} W_j,$$

with

$$W_j = \bigoplus_{w(i_1,\ldots,i_n)=j} V_{i_1} \otimes V_{i_2} \otimes \cdots \otimes V_{i_n},$$

where $w(i_1, \ldots, i_n) = \sharp\{k : i_k = 1\}$. So we have $\dim W_j = \binom{n}{j}(m-1)^j$.

If we define on $X$ the uniform transition probability, that is, $p_u(x, y) = \frac{1}{m}$ for all $x, y \in X$, then the matrix $P_u$ is the matrix $J$ of size $m$.

The eigenvalues of this matrix are 1 (with multiplicity 1) and 0 (with multiplicity $m-1$). The corresponding eigenspaces in $L(X)$ are, respectively, $V_0 \cong \mathbb{C}$ and $V_1 = \{f : X \longrightarrow \mathbb{C} : \sum_{i=0}^{m-1} f(i) = 0\}$.

This means that, by choosing the uniform transition probability on $X$, one gets again the results obtained by considering the Gelfand pair $(S_m, S_{m-1})$.

Also in the case of $X^n$ we can find again the results obtained (see [5]) by considering the Gelfand pair $(S_m \wr S_n, S_{m-1} \wr S_n)$. For $P_u = J$ we have $\lambda_0 = 1$, $\lambda_1 = \cdots = \lambda_{m-1} = 0$. Consider now the transition probability on $X^n$ defined in (4), with $p_1^0 = \cdots = p_n^0 = \frac{1}{n}$. The eigenfunctions $\varphi_I$ associated with the eigenvalue $\frac{1}{n}(n-j)$, with $j = 0, \ldots, n$, are in number of $\binom{n}{j}(m-1)^j$. Moreover

$$\sum_{j=0}^{n} \binom{n}{j}(m-1)^j = \sum_{j=0}^{n} \binom{n}{j}(m-1)^j 1^{n-j} = m^n = \dim L(X^n).$$

For every $j = 0, \ldots, n$, these functions belong to $W_j$ and they are a basis for $W_j$. So $W_j$ is the eigenspace associated with the eigenvalue $\frac{1}{n}(n-j)$.

More generally, consider the case of any reversible transition probability $p$ on $X$. Let $\lambda_0 = 1, \lambda_1, \ldots, \lambda_k$ be the distinct eigenvalues of $P$ and $V_0 \cong \mathbb{C}$, $V_1, \ldots, V_k$ the corresponding eigenspaces. We get

$$L(X^n) \cong (V_0 \oplus V_1 \oplus \cdots \oplus V_k)^{\otimes n}.$$

The eigenfunctions $\varphi_I$ associated with the eigenvalue $\frac{1}{n}\sum_{i=0}^{k} r_i \lambda_i$, with $\sum_{i=0}^{k} r_i = n$, are

$$\binom{r_0 + r_1 + \cdots + r_k}{r_0, \ldots, r_k} \prod_{i=0}^{k} (\dim V_i)^{r_i}$$

and the corresponding eigenspaces are the tensor products in which $r_i$ copies of $V_i$, for $i = 0, 1, \ldots, k$, appear. Moreover, the number of different eigenspaces is equal to the number of integer solutions of the equation

$$r_0 + r_1 + \cdots + r_k = n, \qquad r_i \geq 0,$$

so it is $\binom{k+n}{n}$.

The definition of multinomial coefficient as $\binom{r_0 + r_1 + \cdots + r_k}{r_0, \ldots, r_k} = \frac{(r_0 + \cdots + r_k)!}{r_0! r_1! \cdots r_k!}$ guarantees that

$$\sum_{r_0 + \cdots + r_k = n} \binom{n}{r_0, \ldots, r_k} (\dim V_0)^{r_0} \cdots (\dim V_k)^{r_k} = (\dim V_0 + \cdots + \dim V_k)^n$$

$$= m^n,$$

as we wanted.

## 6. The nested product.

6.1. *General properties.* The nested product of the Markov chains $P_i$'s can be obtained as a particular case of the crested product, by setting $N = \{1, \ldots, n\}$ (so that $C = \varnothing$). The term *nested* comes from the association schemes theory (see [2]).

Consider the product

$$X_1 \times \cdots \times X_n$$

and let $P_i$ be a transition probability on $X_i$. We assume that $p_i$ is in detailed balance with the strict probability measure $\pi_i$, for all $i = 1, \ldots, n$.

The formula of crested product becomes, in this case,

(5) $$P = \sum_{i=1}^{n} p_i^0 (I_1 \otimes \cdots \otimes P_i \otimes J_{i+1} \otimes J_{i+2} \otimes \cdots \otimes J_n).$$

Theorem 4.2 tells us that $P$ is reversible if and only if $P_k$ is symmetric, for every $k > 1$, that is, $\pi_i = \frac{1}{m_i}$ for every $i = 2, \ldots, n$. In this case, $P$ is in detailed balance with the strict probability measure $\pi$ on $X_1 \times \cdots \times X_n$ given by

$$\pi(x_1 \ldots, x_n) = \frac{\pi_1(x_1)}{\prod_{i=2}^{n} m_i}.$$

So let us assume $\pi_i$ to be uniform for every $i = 2, \ldots, n$. The transition probability associated with $P$ is

$$p((x_1, \ldots, x_n), (y_1, \ldots, y_n))$$
$$= \frac{p_1^0 p_1(x_1, y_1)}{m_2 m_3 \cdots m_n} + \sum_{j=2}^{n-1} \frac{\delta((x_1, \ldots, x_{j-1}), (y_1, \ldots, y_{j-1})) p_j^0 p_j(x_j, y_j)}{m_{j+1} \cdots m_n}$$
$$+ \delta((x_1, \ldots, x_{n-1}), (y_1, \ldots, y_{n-1})) p_n^0 p_n(x_n, y_n).$$

As we did in the case of the crossed product, we want to study the spectral theory of the operator $P$ defined in (5).

Let

$$L(X_i) = \bigoplus_{k_i=0}^{r_i} W_{k_i}^i$$

be the spectral decomposition of $L(X_i)$, for all $i = 1, \ldots, n$ and let $\lambda_0^i = 1$, $\lambda_1^i, \ldots, \lambda_{r_i}^i$ be the distinct eigenvalues of $P_i$ associated with these eigenspaces. From Theorem 4.3 we get the following proposition.

PROPOSITION 6.1. *The eigenspaces of $L(X_1 \times \cdots \times X_n)$ are:*

- $L(X_1) \otimes \cdots \otimes L(X_{n-1}) \otimes W_{k_n}^n$, *of eigenvalue $p_n^0 \lambda_{k_n}^n$, for $k_n = 1, \ldots, r_n$, with multiplicity $m_1 \cdots m_{n-1} \dim(W_{k_n}^n)$;*
- $L(X_1) \otimes \cdots \otimes L(X_j) \otimes W_{k_{j+1}}^{j+1} \otimes W_0^{j+2} \otimes \cdots \otimes W_0^n$, *of eigenvalue $p_{j+1}^0 \lambda_{k_{j+1}}^{j+1} + p_{j+2}^0 + \cdots + p_n^0$, with $k_{j+1} = 1, \ldots, r_{j+1}$ and for $j = 1, \ldots, n-2$, with multiplicity $m_1 \cdots m_j \dim(W_{k_{j+1}}^{j+1})$;*
- $W_{k_1}^1 \otimes W_0^2 \otimes \cdots \otimes W_0^n$, *of eigenvalue $p_1^0 \lambda_{k_1}^1 + p_2^0 + \cdots + p_n^0$, for $k_1 = 0, 1, \ldots, r_1$, with multiplicity $\dim(W_{k_1}^1)$.*

Moreover, as in the general case, one can verify that, under the hypothesis that the operators $P_i$ are ergodic, for $i = 1, \ldots, n$, then also the operator $P$ is ergodic.

The application of Proposition 4.5 to the case of the nested product yields the following corollary.

COROLLARY 6.2. *Let $P$ be the nested product of the probabilities $P_i$, with $i = 1, \ldots, n$. Then we have:*

- $U = U_1 \otimes A_2 \otimes \cdots \otimes A_n + \sum_{k=2}^n I_1^{\sigma_1\text{-norm}} \otimes \cdots \otimes I_{k-1}^{\sigma_{k-1}\text{-norm}} \otimes (U_k - A_k) \otimes A_{k+1} \otimes \cdots \otimes A_n;$

- $D = \bigotimes_{i=1}^n D_i$;
- $\Delta = \sum_{i=1}^n p_i^0 (I_1 \otimes \cdots \otimes I_{i-1} \otimes \Delta_i \otimes J_{i+1}^{\text{diag}} \otimes \cdots \otimes J_n^{\text{diag}})$.

6.2. *k-steps transition probability.* The formula that describes the transition probability after $k$ steps in the case of the nested product can be simplified using the base of eigenvectors given in Corollary 6.2 and supposing that the starting point is $0 = (0, \ldots, 0)$.

From the general formula, with the usual notation, we get

$p^{(k)}(0, y)$

$$= \pi(y) \left[ \sum_{z \in X} \left( \sum_{r=2}^n \delta_{\sigma_1}(0, z_1) \cdots \delta_{\sigma_{r-1}}(0, z_{r-1})(u_r - a_r)(0, z_r) \right. \right.$$

$$\times a_{r+1}(0, z_{r+1}) \cdots a_n(0, z_n)$$

$$\left. + u_1(0, z_1) a_2(0, z_2) \cdots a_n(0, z_n) \right) \lambda_z^k$$

$$\times \left( \sum_{r=2}^n \delta_{\sigma_1}(y_1, z_1) \cdots \delta_{\sigma_{r-1}}(y_{r-1}, z_{r-1})(u_r - a_r)(y_r, z_r) \right.$$

$$\times a_{r+1}(y_{r+1}, z_{r+1}) \cdots a_n(y_n, z_n)$$

$$\left. \left. + u_1(y_1, z_1) a_2(y_2, z_2) \cdots a_n(y_n, z_n) \right) \right]$$

$$= \pi(y) \left[ 1 + \sum_{j=2}^n \sum_{\substack{z_j \neq 0 \\ z_i = 0, i \neq j}} \left( \sum_{r=j}^n \frac{1}{\sqrt{\sigma_1(0)} \cdots \sqrt{\sigma_{r-1}(0)}} (u_r - a_r)(0, z_r) \right) \right.$$

$$\times \left( p_r^0 \lambda_{z_r}^r + \sum_{m>r} p_m^0 \right)^k$$

$$\times \left( \sum_{r=j}^n \delta_{\sigma_1}(y_1, 0) \delta_{\sigma_2}(y_2, z_2) \cdots \right.$$

$$\times \delta_{\sigma_{r-1}}(y_{r-1}, z_{r-1})(u_r - a_r)(y_r, z_r)$$

$$\left. \times a_{r+1}(y_{r+1}, z_{r+1}) \cdots a_n(y_n, z_n) \right)$$

$$\left. + \sum_{\substack{z_1 \neq 0 \\ z_i = 0, i > 1}} u_1(0, z_1) \left( p_1^0 \lambda_{z_1}^1 + \sum_{m=2}^n p_m^0 \right)^k u_1(y_1, z_1) \right].$$

Observe that in this case the sum consists of no more than
$$|X_1| + \sum_{i=2}^{n}(|X_i| - 1) = \sum_{i=1}^{n} |X_i| - n + 1$$
nonzero terms.

EXAMPLE 6.3. We want to express the $k$th step transition probability in the case $n = 2$. So consider the product $X \times Y$, with $X = \{0, 1, \ldots, m\}$ and $Y = \{0, 1, \ldots, n\}$. Let
$$L(X) = \bigoplus_{j=0}^{r} V_j \quad \text{and} \quad L(Y) = \bigoplus_{i=0}^{s} W_i$$
be the spectral decomposition of the spaces $L(X)$ and $L(Y)$, respectively. Let $\lambda_0 = 1, \lambda_1, \ldots, \lambda_r$ and $\mu_0 = 1, \mu_1, \ldots, \mu_s$ be the distinct eigenvalues of $P_X$ and $P_Y$, respectively. Then the eigenspaces of $L(X \times Y)$ are $L(X) \otimes W_i$, for $i = 1, \ldots, s$, with dimension $(m+1) \dim W_i$ and associated eigenvalue $p_Y^0 \mu_i$, and $V_j \otimes W_0$, for $j = 0, \ldots, r$, with dimension $\dim V_j$ and associated eigenvalue $p_X^0 \lambda_j + p_Y^0$.

The expression of $U$ becomes
$$U = I_X^{\sigma_X\text{-norm}} \otimes (U_Y - A_Y) + U_X \otimes A_Y.$$
In particular, let $\{v^0, v_1^1, \ldots, v_{\dim(V_1)}^1, \ldots, v_1^r, \ldots, v_{\dim(V_r)}^r\}$ and $\{w^0, w_1^1, \ldots, w_{\dim(W_1)}^1, \ldots, w_1^s, \ldots, w_{\dim(W_s)}^s\}$ be the eigenvectors of $P_X$ and $P_Y$, respectively, that is, they represent the columns of the matrices $U_X$ and $U_Y$.

Then, the columns of the matrix $U$ corresponding to the elements $(i, 0) \in \{0, \ldots, m\} \times \{0, \ldots, n\}$ are the eigenvectors $v^i \otimes (1, \ldots, 1)$ with eigenvalue $p_X^0 \lambda_i + p_Y^0$. On the other hand, the columns corresponding to the elements $(i, j) \in \{0, \ldots, m\} \times \{0, \ldots, n\}$, with $j = 1, \ldots, n$, are the eigenvectors
$$\Big(0, \ldots, 0, \underbrace{\frac{1}{\sqrt{\sigma_X(i)}}}_{i\text{th place}}, 0, \ldots, 0\Big) \otimes w^j$$
whose eigenvalue is $p_Y^0 \mu_j$. As a consequence, only $m+1+n$ of these eigenvectors can be nonzero in the first coordinate, so the probability $p^{(k)}((0,0),(x,y))$ can be expressed as a sum of $m+1+n$ nonzero terms; moreover, these terms become $m+1$ if $x \neq 0$. We have

$p^{(k)}((0,0),(x,y))$
$$= \pi((x,y))\Bigg(\sum_{i=0}^{m} v^i(0) v^i(x)(p_X^0 \lambda_i + p_Y^0)^k$$
$$+ \frac{1}{\sqrt{\sigma_X(0)\sigma_X(x)}} \sum_{j=1}^{n} w^j(0)\delta_0(x) w^j(y)(p_Y^0 \mu_j)^k\Bigg)$$

$$= \frac{\sigma_X(x)}{n+1}\left[\sum_{i=0}^{r}\left(\sum_{a=1}^{\dim(V_i)} v_a^i(0)v_a^i(x)\right)(p_X^0\lambda_i + p_Y^0)^k\right.$$
$$\left.+ \sum_{j=1}^{s}\left(\frac{1}{\sqrt{\sigma_X(0)\sigma_X(x)}}\sum_{b=1}^{\dim(W_j)} w_b^j(0)\delta_0(x)w_b^j(y)\right)(p_Y^0\mu_j)^k\right].$$

6.3. *The insect.* It is clear that the product $X_1 \times \cdots \times X_n$ can be regarded as the rooted tree $T$ of depth $n$, such that the root has degree $m_1$, each vertex of the first level has $m_2$ children and in general each vertex of the $i$th level of the tree has $m_{i+1}$ children, for every $i = 1, \ldots, n-1$. We denote the $i$th level of the tree by $L_i$. In this way, every vertex $x \in L_i$ can be regarded as a word $x = x_1 \cdots x_i$, where $x_j \in \{0, 1, \ldots, m_j - 1\}$.

We want to show that the nested product of Markov chains is the generalization of the "insect problem" studied by Figà-Talamanca in [11].

Let us imagine that an insect lives in a leaf $x \in L_n$ and that it performs a simple random walk on the graph $T$ starting from $x$.

Then there exists a probability distribution $\mu_x$ on $L_n$ such that, for every $y \in L_n$, $\mu_x(y)$ is the probability that $y$ is the first point in $L_n$ visited by the insect in the random walk. If we put $\overline{p}(x,y) = \mu_x(y)$, then we get a stochastic matrix $\overline{P} = (\overline{p}(x,y))_{x,y \in L_n}$. Since the random walk is $Aut(T)$-invariant, we can suppose that the random walk starts at the leftmost vertex, that we will call $x_0 = (0, \ldots, 0)$. We recall that $Aut(T)$ is the group of all automorphisms of $T$, given by the iterated wreath product $S_{m_n} \wr S_{m_{n-1}} \wr \cdots \wr S_{m_1}$. We want to study this Markov chain defined on $L_n$.

Set $\xi_n = \varnothing$ and $\xi_i = 00\ldots 0$ ($n-i$ times). For $j \geq 0$, let $\alpha_j$ be the probability that the insect reaches $\xi_{j+1}$ given that $\xi_j$ is reached at least once. This definition implies $\alpha_0 = 1$ and $\alpha_1 = \frac{1}{m_n+1}$. In fact, with probability 1, the insect reaches the vertex $\xi_1$ at the first step and, starting from $\xi_1$, with probability $\frac{1}{m_n+1}$ it reaches $\xi_2$, while with probability $\frac{m_n}{m_n+1}$ it returns to $L_n$. Finally, we have $\alpha_n = 0$. For $1 < j < n$, there is the following recursive relation:

$$\alpha_j = \frac{1}{m_{n+1-j}+1} + \alpha_{j-1}\alpha_j\frac{m_{n+1-j}}{m_{n+1-j}+1}.$$

In fact, starting at $\xi_j$, with probability $\frac{1}{m_{n+1-j}+1}$ the insect reaches in one step $\xi_{j+1}$, otherwise with probability $\frac{m_{n+1-j}}{m_{n+1-j}+1}$ it reaches $\xi_{j-1}$ or one of its brothers; then, with probability $\alpha_{j-1}$ it reaches again $\xi_j$ and one starts the recursive argument.

The solution, for $1 \leq j \leq n-1$, is given by

$$\alpha_j = \frac{1 + m_n + m_nm_{n-1} + m_nm_{n-1}m_{n-2} + \cdots + m_nm_{n-1}m_{n-2}\cdots m_{n-j+2}}{1 + m_n + m_nm_{n-1} + m_nm_{n-1}m_{n-2} + \cdots + m_nm_{n-1}m_{n-2}\cdots m_{n-j+1}}$$

$$= 1 - \frac{m_n m_{n-1} m_{n-2} \cdots m_{n-j+1}}{1 + m_n + m_n m_{n-1} + m_n m_{n-1} m_{n-2} + \cdots + m_n m_{n-1} m_{n-2} \cdots m_{n-j+1}}.$$

Moreover, we have

$$\overline{p}(x_0, x_0) = \frac{1}{m_n}(1 - \alpha_1) + \frac{1}{m_n m_{n-1}} \alpha_1 (1 - \alpha_2) + \cdots$$

$$+ \frac{1}{m_n m_{n-1} \cdots m_2} \alpha_1 \alpha_2 \cdots \alpha_{n-2}(1 - \alpha_{n-1}) + \frac{1}{m_n \cdots m_1} \alpha_1 \cdots \alpha_{n-1}.$$

Indeed the $j$th summand is the probability of returning back to $x_0$ if the corresponding random walk in $T$ reaches $\xi_j$ but not $\xi_{j+1}$. It is not difficult to compute $\overline{p}(x_0, x)$, where $x$ is a point at distance $j$ from $x_0$. For $j = 1$, we clearly have $\overline{p}(x_0, x_0) = \overline{p}(x_0, x)$. We observe that, for $j > 1$, to reach $x$ one is forced to first reach $\xi_j$, so that we have

$$\overline{p}(x_0, x) = \frac{1}{m_n \cdots m_{n-j+1}} \alpha_1 \alpha_2 \cdots \alpha_{j-1}(1 - \alpha_j) + \cdots$$

$$+ \frac{1}{m_n \cdots m_2} \alpha_1 \alpha_2 \cdots \alpha_{n-2}(1 - \alpha_{n-1}) + \frac{1}{m_n \cdots m_1} \alpha_1 \alpha_2 \cdots \alpha_{n-1}.$$

Since the random walk is invariant with respect to the action of $Aut(T)$, which acts isometrically on the tree, we get the same formula for any pair of vertices $x, y \in L_n$ such that $d(x, y) = j$.

PROPOSITION 6.4. *The stochastic matrix*

$$p((x_1, \ldots, x_n), (y_1, \ldots, y_n))$$

$$= \frac{p_1^0 p_1(x_1, y_1)}{m_2 m_3 \cdots m_n} + \sum_{j=2}^{n-1} \frac{\delta((x_1, \ldots, x_{j-1}), (y_1, \ldots, y_{j-1})) p_j^0 p_j(x_j, y_j)}{m_{j+1} \cdots m_n}$$

$$+ \delta((x_1, \ldots, x_{n-1}), (y_1, \ldots, y_{n-1})) p_n^0 p_n(x_n, y_n),$$

*defined in (5), gives rise to the Insect Markov chain on $L_n$, regarded as $X_1 \times \cdots \times X_n$, choosing $p_i^0 = \alpha_1 \alpha_2 \cdots \alpha_{n-i}(1 - \alpha_{n-i+1})$ for $i = 1, \ldots, n-1$ and $p_n^0 = 1 - \alpha_1$ and the transitions probabilities $p_j$'s to be uniform for all $j = 1, \ldots, n$.*

PROOF. Set, for every $i = 1, \ldots, n-1$,

$$p_i^0 = \alpha_1 \alpha_2 \cdots \alpha_{n-i}(1 - \alpha_{n-i+1})$$

and $p_n^0 = 1 - \alpha_1$. Moreover, assume that the probability $p_i$ on $X_i$ is uniform, that is,

$$P_i = J_i.$$

If $d(x_0, x) = n$, then we get
$$p(x_0, x) = \frac{\alpha_1 \alpha_2 \cdots \alpha_{n-1}}{m_1 m_2 \cdots m_n}.$$
If $d(x_0, x) = j > 1$, that is, $x_i^0 = x_i$ for all $i = 1, \ldots, n-j$, then
$$p(x_0, x) = \frac{\alpha_1 \alpha_2 \cdots \alpha_{n-1}}{m_1 m_2 \cdots m_n} + \sum_{i=1}^{n-j} \frac{\alpha_1 \cdots \alpha_{n-i-1}(1-\alpha_{n-i})}{m_n \cdots m_{i+2} m_{i+1}}.$$
Finally, if $x = x_0$, we get
$$p(x_0, x_0) = \frac{\alpha_1 \alpha_2 \cdots \alpha_{n-1}}{m_1 m_2 \cdots m_n} + \sum_{i=1}^{n-2} \frac{\alpha_1 \cdots \alpha_{n-i-1}(1-\alpha_{n-i})}{m_n \cdots m_{i+2} m_{i+1}} + \frac{(1-\alpha_1)}{m_n}.$$
This completes the proof. □

The decomposition of the space $L(L_n) = L(X_1 \times \cdots \times X_n)$ under the action of $Aut(T)$ is known (see [6]). Denote $Z_0 \cong \mathbb{C}$ the trivial representation and, for every $j = 1, \ldots, n$, define the following subspace:
$$Z_j = \left\{ f \in L(L_n) : f = f(x_1, \ldots, x_j), \sum_{i=0}^{m_j - 1} f(x_1, \ldots, x_{j-1}, i) \equiv 0 \right\}$$
of dimension $m_1 \cdots m_{j-1}(m_j - 1)$. In virtue of the correspondence between $Aut(T)$-invariant operators and bi-$Stab_{Aut(T)}(0^n)$-invariant functions, the corresponding eigenvalues are given by the spherical Fourier transform of the convolver that represents $\overline{P}$, namely
$$\lambda_j = \sum_{x \in L_n} \overline{P}(x_0, x) \phi_j(x),$$
where $\phi_j$ is the $j$th spherical function, for all $j = 0, 1, \ldots, n$. It is easy to verify that one gets:

- $\lambda_0 = 1$;
- $\lambda_j = 1 - \alpha_1 \alpha_2 \cdots \alpha_{n-j}$, for every $j = 1, \ldots, n-1$;
- $\lambda_n = 0$.

In particular, if we set
$$p_i^0 = \alpha_1 \alpha_2 \cdots \alpha_{n-i}(1 - \alpha_{n-i+1})$$
for every $i = 1, \ldots, n-1$, with $p_n^0 = 1 - \alpha_1$ and $P_i = J_i$ for every $i = 1, \ldots, n$, the eigenspaces given for $L(X_1 \times \cdots \times X_n)$ in Proposition 6.1 are exactly the $Z_j$'s with the corresponding eigenvalues.

Let us prove that the eigenvalues that we have obtained in Proposition 6.1 coincide with the eigenvalues corresponding to the eigenspaces $Z_0, Z_1, \ldots, Z_n$.

We want to get these eigenvalues by using the formulas given in Proposition 6.1 for the eigenvalues of the nested product $P$ by setting $P_i = J_i$, then $p_i^0 = \alpha_1 \alpha_2 \cdots \alpha_{n-i}(1 - \alpha_{n-i+1})$ for $i = 1, \ldots, n-1$ and $p_n^0 = 1 - \alpha_1$. First of all, we observe that the eigenvalues of the operator $P_i$ are 1, with multiplicity 1 and 0, with multiplicity $m_i - 1$. So we get $L(X_i) = W_0^i \oplus W_1^i$, with $\dim(W_1^i) = m_i - 1$, for all $i = 1, \ldots, n$. Following the formulas that we have given, the eigenspaces of $P$ are:

- $L(X_1) \otimes L(X_2) \otimes \cdots \otimes L(X_{n-1}) \otimes W_1^n$;
- $L(X_1) \otimes L(X_2) \otimes \cdots \otimes L(X_{n-j-1}) \otimes W_1^{n-j} \otimes W_0^{n-j+1} \otimes \cdots \otimes W_0^n$, for every $j = 1, \ldots, n-1$;
- $W_0^1 \otimes W_0^2 \otimes \cdots \otimes W_0^n$.

The corresponding eigenvalues are:

- $p_n^0 \lambda_1^n = 0$;
- $\sum_{i=n-j+1}^{n} p_i^0$, for every $j = 1, \ldots, n-1$;
- $\sum_{i=1}^{n} p_i^0 = 1$.

We need to prove that, for every $j = 1, \ldots, n-1$, the eigenvalue $\sum_{i=n-j+1}^{n} p_i^0$ is equal to $1 - \alpha_1 \alpha_2 \cdots \alpha_j$. We prove the assertion by induction on $j$.

If $j = 1$, we have $p_n^0 = 1 - \alpha_1$. Now suppose the assertion to be true for $j$ and show that it holds also for $j + 1$. We get

$$\sum_{i=n-j}^{n} p_i^0 = \sum_{i=n-j+1}^{n} p_i^0 + p_{n-j}^0 = 1 - \alpha_1 \alpha_2 \cdots \alpha_j + \alpha_1 \cdots \alpha_j (1 - \alpha_{j+1})$$
$$= 1 - \alpha_1 \cdots \alpha_j \alpha_{j+1}.$$

**7. The second crested product.** In this section we define a different kind of product of two spaces $X$ and $Y$, that we will call the second crested product. In fact it contains, as particular cases, the crossed product and the nested product described in Sections 5 and 6, respectively. We will study a Markov chain $P$ on the set $\Theta_k$ of functions from $X$ to $Y$ whose domains are $k$-subsets of $X$, giving the spectrum and the relative eigenspaces.

7.1. *General theory.* Let $X$ be a finite set of cardinality $n$, say $X = \{1, 2, \ldots, n\}$. For every $k = 1, \ldots, n$, denote by $\Omega_k$ the set of $k$-subsets of $X$, so that $|\Omega_k| = \binom{n}{k}$.

Now let $Y$ be a finite set and let $Q$ be a transition matrix on $Y$, which is in detailed balance with the strict probability $\tau$. Let $\lambda_0 = 1, \lambda_1, \ldots, \lambda_m$ be the distinct eigenvalues of $Q$ and denote by $W_j$ the corresponding eigenspaces, for every $j = 0, 1, \ldots, m$, so that the following spectral decomposition holds:

$$L(Y) = \bigoplus_{j=0}^{m} W_j.$$

Moreover, assume that the dimension of the eigenspace associated with the eigenvalue 1 is 1 and set $\dim(W_j) = m_j$, for every $j = 1, \ldots, m$.

Recall that the eigenspace $W_0$ is generated by the vector

$$\underbrace{(1, \ldots, 1)}_{|Y| \text{ times}}$$

and that $W_j$ is orthogonal to $W_0$ with respect to the scalar product $\langle \cdot, \cdot \rangle_\tau$, for every $j = 1, \ldots, m$.

For every $k = 1, \ldots, n$, consider the space

$$\Theta_k = \{(A, \theta) : A \in \Omega_k \text{ and } \theta \in Y^A\},$$

that is, the space of functions $\theta$ whose domain $A = dom(\theta)$ is a $k$-subset of $X$ and which take values in $Y$.

The set $\Theta = \coprod_{k=0}^n \Theta_k$ is a poset with respect to the relation $\subseteq$ defined in the following way:

$$\varphi \subseteq \chi \text{ if } dom(\varphi) \subseteq dom(\chi) \text{ and } \varphi = \chi|_{dom(\varphi)}.$$

The Markov chain $P$ on $\Theta_k$ that we are going to define can be regarded as follows. Let $0 < p_0 < 1$ a real number. Then, starting from a function $\theta \in \Theta_k$, with probability $p_0$ we can reach a function $\varphi \in \Theta_k$ having the same domain as $\theta$ and that can differ from $\theta$ at most in one image, according with the probability $Q$ on $Y$.

On the other hand, with probability $1 - p_0$ we can reach in one step a function $\varphi \in \Theta_k$ whose domain intersects the domain of $\theta$ in $k-1$ elements (on which the functions coincide), and in such a way that the image of the $k$th element of the domain of $\varphi$ is uniformly chosen.

Note that $P$ defines a Markovian operator on the space $L(\Theta_k)$ of all complex functions defined on $\Theta_k$.

When $Y$ is the ultrametric space, the Markov chain $P$ represents the so-called multi-insect, which generalizes the Insect Markov chain already studied. In particular if $|X| = n$, we consider $k$ insects living in $k$ different subtrees and moving only one per each step in such a way that their distance is preserved, giving rise to a Markov chain on the space of all possible configurations of $k$ insects having this property.

In fact each element in $\Theta_k$ can be regarded as a configuration of $k$ insects and vice versa. For example, let $\theta \in \Theta_k$ be a function such that $dom(\theta) = \{x_1, \ldots, x_k\}$ and $\theta(x_i) = y_i$, with $x_i \in X$ and $y_i \in Y$ for all $i = 1, \ldots, k$. Then the corresponding configuration of $k$ insects has an insect at each leaf $(x_i, y_i)$. They live in all different subtrees since $x_i \neq x_j$ for $i \neq j$.

We observe that the cardinality of this space is $\binom{n}{k}|Y|^k$. This space can be regarded as the variety of subtrees (see [4]) of branch indices $(k, 1)$ in the rooted tree $(n, |Y|)$.

If $\theta, \varphi \in \Theta_k$, with domains $A$ and $B$ respectively, then define the matrix $\Delta$, indexed by $\Theta_k$, whose entries are

$$\Delta_{\theta,\varphi} = \begin{cases} 1, & \text{if } |A \cap B| = k-1 \text{ and } \theta|_{A \cap B} = \varphi|_{A \cap B}, \\ 0, & \text{otherwise.} \end{cases}$$

Observe that the matrix $\Delta$ is symmetric.

The operator $P$ can be expressed in terms of the operator associated with $\Delta$ and of another operator $M$ as

$$(6) \qquad P = p_0 M + (1 - p_0) \frac{\Delta}{norm(\Delta)},$$

where $M$ describes the situation in which the domain is not changed and only one of the images of the function $\theta \in \Theta_k$ is changed according with the probability $Q$ on $Y$. An analytic expression for $M$ will be presented below. On the other hand, $\Delta$ describes the situation in which we pass from a function whose domain is $A$ to a function whose domain is $A \sqcup \{i\} \setminus \{j\}$, with $i \notin A$ and $j \in A$, and we choose uniformly the image in $Y$ of the element $i$. So the action of $\Delta$ on $\Omega_k$ is an analogue of the Bernoulli–Laplace diffusion model. By $norm(\Delta)$ we indicate the number of nonzero entries in each row of the matrix associated with $\Delta$.

It is easy to check that $M$ is in detailed balance with the strict probability measure defined as

$$\tau_M(\theta) = \frac{1}{\binom{n}{k}} \prod_{i \in A} \tau(\theta(i)),$$

where $\theta \in \Theta_k$ and $dom(\theta) = A$. On the other hand, it follows from the definition of the Markov chain $\Delta$ that the weighted graph associated with $\Delta$ is connected. From this and from the fact that the nonzero entries of $\Delta$ are all equal to 1, we can deduce that $\Delta$ is reversible and in detailed balance with a uniform probability measure. This forces $\tau_M$ to be uniform and so we have to assume that $\tau$ is uniform on $Y$ and the matrix $Q$ is symmetric.

In this way, $P$ is in detailed balance with the uniform probability measure $\pi$ such that $\pi(\theta) = \frac{1}{\binom{n}{k}|Y|^k}$, for every $\theta \in \Theta_k$. This choice of $\tau$ guarantees that, if $f$ is any function in $W_j$, with $j = 1, \ldots, m$, then $\sum_{y \in Y} f(y) = 0$.

The spectral theory of the operator $M$ has been studied in Section 5. In fact, it corresponds to choosing, with probability $\frac{1}{k}$, only one element of the domain and to changing the corresponding image with respect to the probability $Q$ on $Y$, fixing the remaining ones. So we focus our attention to investigate the spectral theory of the operator $\Delta$.

Let us introduce two differential operators.

DEFINITION 7.1. (1) For every $k = 2, \ldots, n$ the operator $D_k : L(\Theta_k) \longrightarrow L(\Theta_{k-1})$ is defined by

$$(D_k F)(\varphi) = \sum_{\theta \in \Theta_k : \varphi \subseteq \theta} F(\theta),$$

for every $F \in L(\Theta_k)$ and $\varphi \in \Theta_{k-1}$.

(2) For $k = 1, \ldots, n$ the operator $D_k^* : L(\Theta_{k-1}) \longrightarrow L(\Theta_k)$ is defined by

$$(D_k^* F)(\theta) = \sum_{\varphi \in \Theta_{k-1} : \varphi \subseteq \theta} F(\varphi),$$

for every $F \in L(\Theta_{k-1})$ and $\theta \in \Theta_k$.

Observe that the operator $D_k^*$ is the adjoint of $D_k$. The following decomposition holds:

$$L(\Theta_k) = L\left(\coprod_{A \in \Omega_k} Y^A\right) = \bigoplus_{A \in \Omega_k} L(Y^A).$$

In order to get a basis for the space $L(Y^A)$, for every $A \in \Omega_k$, we introduce some special functions that we will call fundamental functions.

DEFINITION 7.2. Suppose that $A \in \Omega_k$ and that $F^j \in L(Y)$ for every $j \in A$. Suppose also that each $F^j$ belongs to an eigenspace of $Q$ and set $a_i = |\{j \in A : F^j \in W_i\}|$. Then the tensor product $\bigotimes_{j \in A} F^j$ will be called a fundamental function of type $\underline{a} = (a_0, a_1, \ldots, a_m)$ in $L(Y^A)$.

In other words, we have

$$\left(\bigotimes_{j \in A} F^j\right)(\theta) = \prod_{j \in A} F^j(\theta(j)),$$

for every $\theta \in Y^A$. We also set $\ell(\underline{a}) = a_1 + \cdots + a_m = k - a_0$.

The introduction of the fundamental functions allows to give a useful expression for the operators $M$ and $\Delta$.

If $F \in L(Y^A) \subseteq L(\Theta_k)$ is the fundamental function $F = \bigotimes_{j \in A} F^j$, with $|A| = k$ and $F^j : Y \longrightarrow \mathbb{C}$, then $MF = \frac{1}{k} \sum_{j \in A}[(\bigotimes_{i \in A, i \neq j} F^i) \otimes QF^j]$. So, if $\theta \in \Theta_k$ and $dom(\theta) = A$, we get

$$(MF)(\theta) = \frac{1}{k} \sum_{j \in A} \left[\prod_{i \in A, i \neq j} F^i(\theta(i)) \left(\sum_{y \in Y} q(\theta(j), y) F^j(y)\right)\right].$$

Analogously one has $(\Delta F)(\theta) = \sum_\varphi F(\varphi)$, where the sum is over all $\varphi \in \Theta_k$ such that $dom(\varphi) \cap dom(\theta) = k - 1$ and $\varphi \equiv \theta$ on $dom(\varphi) \cap dom(\theta)$. If

$A = (dom(\theta) \cap A) \sqcup \{i\}$ (we denote by $\sqcup$ the disjoint union), then

$$\left(\Delta\left(\bigotimes_{j \in A} F^j\right)\right)(\theta) = \sum_{\varphi} \bigotimes_{j \in A} F^j(\varphi) = \prod_{j \in dom(\varphi) \cap A} F^j(\theta(j)) \left(\sum_{y \in Y} F^i(y)\right).$$

Denote $P_{k,\underline{a},A}$ the subspace of $L(Y^A)$ spanned by the fundamental functions of type $\underline{a}$ and

$$P_{k,\underline{a}} = \bigoplus_{A \in \Omega_k} P_{k,\underline{a},A}.$$

LEMMA 7.3. $D_k$ maps $P_{k,\underline{a}}$ to $P_{k-1,\underline{a}'}$, where $\underline{a}' = (a_0 - 1, a_1, \ldots, a_m)$. Conversely $D_k^*$ maps $P_{k-1,\underline{a}'}$ to $P_{k,\underline{a}}$.

PROOF. Let $F$ be a fundamental function of type $\underline{a}$ in $L(Y^A)$ and let $B \subset A$ such that $A = B \sqcup \{i\}$. Then for every $\varphi \in Y^B$, we have

$$(D_k F)(\varphi) = \sum_{\theta \in Y^A : \varphi \subseteq \theta} F(\theta)$$

$$= \sum_{\theta \in Y^A : \varphi \subseteq \theta} \prod_{j \in A} F^j(\theta(j))$$

$$= \left(\sum_{y \in Y} F^i(y)\right) \prod_{j \in B} F^j(\varphi(j)).$$

The value of $\sum_{y \in Y} F^i(y)$ is zero if $F^i \in W_j$ for $j = 1, \ldots, m$ and so $D_k F \equiv 0$ if $a_0 = 0$. If $F^i \in W_0$, then $D_k F \in P_{k-1,\underline{a}'}$.

Analogously, let $F \in P_{k-1,\underline{a}',B}$ with $B \in \Omega_{k-1}$. Then for every $\theta \in Y^A$, $A = B \sqcup \{i\}$, one has

$$(D_k^* F)(\theta) = \sum_{\varphi \in Y^B : \varphi \subseteq \theta} F(\varphi)$$

$$= \prod_{j \in B} F^j(\varphi(j))$$

$$= F^i(\theta(i)) \prod_{j \in B} F^j(\theta(j)),$$

where $F^i \equiv 1$ on $Y$ (and so $F^i \in W_0$). □

The restriction of $D_k$ to $P_{k,\underline{a}}$ will be denoted by $D_{k,\underline{a}}$ and the restriction of $D_k^*$ to $P_{k-1,\underline{a}'}$ by $D_{k,\underline{a}}^*$.

The study of the compositions of the operators $D_{k,\underline{a}}$ and $D^*_{k,\underline{a}}$ plays a central role. In fact it will be shown that the eigenspaces of these operators are also eigenspaces for $\Delta$. Consider, for example, $D_{k+1}D^*_{k+1}$ applied to a function $F \in L(\Theta_k)$ and calculated on $\theta \in \Theta_k$. The functions $\varphi \in \Theta_{k+1}$ such that $\varphi \supseteq \theta$ are in number of $|Y|(n-k)$. Each of them covers $k+1$ functions in $\Theta_k$; one of them is the function $\theta$, the other ones are functions in $\Theta_k$ whose domains differ by the domain of $\theta$ of an element and coincide on their intersection. These functions are in number of $|Y|(n-k)k$ and they correspond to functions that one can reach starting from $\theta$ in the Markov chain described by $\Delta$. From this it follows that $norm(\Delta) = |Y|(n-k)k$.

LEMMA 7.4. *Let $F \in P_{k,\underline{a},A}$, with $A \in \Omega_k$. Then*

$$D^*_{k,\underline{a}}D_{k,\underline{a}} = |Y|(k - \ell(\underline{a}))I + Q_{k,\underline{a}},$$

*where $Q_{k,\underline{a}}$ is defined by setting*

(7) $$(Q_{k,\underline{a}}F)(\theta) = \begin{cases} 0, & \text{if } F^i \notin W_0, \\ |Y|F(\overline{\theta}), & \text{if } F^i \in W_0, \end{cases}$$

*for every $\theta \in \Theta_k$ such that $|dom(\theta) \cap A| = k-1$ and $A \setminus dom(\theta) = \{i\}$. We denote by $\overline{\theta}$ the function in $\Theta_k$ whose domain is $A$ and such that $\overline{\theta}|_{A\setminus\{i\}} = \theta$ and $\overline{\theta}(i) = \theta(i_0)$, where $dom(\theta) \setminus A = \{i_0\}$.*

PROOF. Take $F \in P_{k,\underline{a},A}$ and $\theta \in \Theta_k$. We have

$$(D^*_{k,\underline{a}}D_{k,\underline{a}}F)(\theta) = \sum_{\varphi \in \Theta_{k-1}: \varphi \subseteq \theta} (D_{k,\underline{a}}F)(\varphi)$$
$$= \sum_{\varphi \in \Theta_{k-1}: \varphi \subseteq \theta} \sum_{\substack{\omega \in \Theta_k: \omega \supseteq \varphi, \\ dom(\omega) = A}} F(\omega).$$

If $dom(\theta) = A$, then we get

$$(D^*_{k,\underline{a}}D_{k,\underline{a}}F)(\theta) = \sum_{j \in A}\left(\sum_{y \in Y} F^j(y)\right) \prod_{t \in A\setminus\{j\}} F^t(\theta(t))$$
$$= |Y|(k - \ell(\underline{a})) \prod_{t \in A} F^t(\theta(t))$$
$$= |Y|(k - \ell(\underline{a}))F(\theta),$$

where the second equality follows from the fact that $\sum_{y \in Y} F^j(y) = |Y|$ if $F^j \in W_0$ and $\sum_{y \in Y} F^j(y) = 0$ whenever $F^j \notin W_0$.

On the other hand, if $|dom(\theta) \cap A| = k-1$, with $A \setminus dom(\theta) = \{i\}$, then

$$(D_{k,\underline{a}}^* D_{k,\underline{a}} F)(\theta) = \left(\sum_{y \in Y} F^i(y)\right) \prod_{j \in A \setminus \{i\}} F^j(\theta(j))$$

$$= \begin{cases} 0, & \text{if } F^i \notin W_0, \\ |Y|F(\overline{\theta}), & \text{if } F^i \in W_0, \end{cases}$$

which is just the definition of $Q_{k,\underline{a}}$. □

LEMMA 7.5. *Let $F \in P_{k,\underline{a}',A}$, with $A \in \Omega_k$. Then*

$$D_{k+1,\underline{a}} D_{k+1,\underline{a}}^* = |Y|(n-k)I + Q_{k,\underline{a}},$$

*where $Q_{k,\underline{a}}$ is defined as in (7).*

PROOF. Take $F \in P_{k,\underline{a}',A}$ and $\theta \in \Theta_k$. We have

$$(D_{k+1,\underline{a}} D_{k+1,\underline{a}}^* F)(\theta) = \sum_{\varphi \in \Theta_{k+1}: \theta \subseteq \varphi} (D_{k+1,\underline{a}}^* F)(\varphi)$$

$$= \sum_{\varphi \in \Theta_{k+1}: \theta \subseteq \varphi} \sum_{\substack{\omega \in \Theta_k: \omega \supseteq \varphi, \\ dom(\omega) = A}} F(\omega).$$

If $dom(\theta) = A$, then we get

$$(D_{k+1,\underline{a}} D_{k+1,\underline{a}}^* F)(\theta) = \sum_{j \in A^C} \sum_{y \in Y} F(\theta)$$

$$= |Y|(n-k)F(\theta).$$

On the other hand, if $|dom(\theta) \cap A| = k-1$, with $A \setminus dom(\theta) = \{i\}$, then

$$(D_{k+1,\underline{a}} D_{k+1,\underline{a}}^* F)(\theta) = \left(\sum_{y \in Y} F^i(y)\right) \prod_{j \in A \setminus \{i\}} F^j(\theta(j))$$

$$= \begin{cases} 0, & \text{if } F^i \notin W_0, \\ |Y|F(\overline{\theta}), & \text{if } F^i \in W_0, \end{cases}$$

$$= (Q_{k,\underline{a}} F)(\theta).$$

This completes the proof. □

The following corollary easily follows.

COROLLARY 7.6. *Let $F \in P_{k,\underline{a}',A}$, with $A \in \Omega_k$. Then*

$$D_{k+1,\underline{a}} D_{k+1,\underline{a}}^* - D_{k,\underline{a}'}^* D_{k,\underline{a}'} = |Y|(n + \ell(\underline{a}) - 2k)I.$$

Consider now the operator $D_{k,\underline{a}}\colon P_{k,\underline{a}} \longrightarrow P_{k-1,\underline{a}'}$.

DEFINITION 7.7. For $0 \leq \ell(\underline{a}) \leq k \leq n$, set
$$P_{k,\underline{a},k} = Ker(D_{k,\underline{a}})$$
and inductively, for $k \leq h \leq n$, set
$$P_{h,\underline{a},k} = D^*_{h,\underline{a}} P_{h-1,\underline{a}',k}.$$

These spaces have a fundamental importance because they exactly constitute the eigenspaces of the operator $\Delta$. This will be a consequence of the following proposition.

PROPOSITION 7.8. *$P_{h,\underline{a}',k}$ is an eigenspace for the operator $D_{h+1,\underline{a}} D^*_{h+1,\underline{a}}$ and the corresponding eigenvalue is $|Y|(n+\ell(\underline{a})-k-h)(h-k+1)$.*

PROOF. We prove the assertion by induction on $h$. If $h=k$, from the last corollary we get $D_{k+1,\underline{a}} D^*_{k+1,\underline{a}}|_{P_{k,\underline{a}',k}} = |Y|(n+\ell(\underline{a})-2k)I$, since $D_{k,\underline{a}'} P_{k,\underline{a}',k} = 0$ by definition of $P_{k,\underline{a}',k}$.

Now suppose the lemma to be true for $k \leq t \leq h$ and recall that, by definition, we have $P_{h+1,\underline{a}',k} = D^*_{h+1,\underline{a}'} P_{h,\underline{a}'',k}$. Moreover, Corollary 7.6 gives
$$D_{h+2,\underline{a}} D^*_{h+2,\underline{a}} - D^*_{h+1,\underline{a}'} D_{h+1,\underline{a}'} = |Y|(n+\ell(\underline{a})-2(h+1))I.$$
So we get
$$\begin{aligned} D_{h+2,\underline{a}} D^*_{h+2,\underline{a}}|_{P_{h+1,\underline{a}',k}} &= D^*_{h+1,\underline{a}'}|_{D_{h+1,\underline{a}'} D^*_{h+1,\underline{a}'} P_{h,\underline{a}'',k}} \\ &\quad + |Y|(n+\ell(\underline{a})-2(h+1)) P_{h+1,\underline{a}',k} \\ &= |Y|(n+\ell(\underline{a})-k-h)(h-k+1) D^*_{h+1,\underline{a}'} P_{h,\underline{a}'',k} \\ &\quad + |Y|(n+\ell(\underline{a})-2(h+1)) P_{h+1,\underline{a}',k} \\ &= |Y|(n+\ell(\underline{a})-k-h-1)(h-k+2) P_{h+1,\underline{a}',k}, \end{aligned}$$
where the second equality follows from the inductive hypothesis and the third one from an easy computation. This completes the proof. □

COROLLARY 7.9. *$P_{h,\underline{a}',k}$ is an eigenspace for $\Delta$ of eigenvalue $|Y|(n+\ell(\underline{a})-k-h)(h-k+1) - |Y|(n-h)$.*

PROOF. It suffices to observe that the operator $Q_{h,\underline{a}}$ defined in (7) coincides with the operator $\Delta$ on the space $P_{h,\underline{a}}$ and then the assertion follows from Lemma 7.5 and Proposition 7.8. □

In particular, after normalizing the matrix $\Delta$ we obtain $\frac{\Delta}{norm(\Delta)}$ and the corresponding eigenvalue is $\frac{1}{|Y|(n-h)h}(|Y|(n+\ell(\underline{a})-k-h)(h-k+1) - |Y|(n-h))$.

The following lemma holds.

LEMMA 7.10. *Given $\ell(\underline{a})$ and $h$, then, for $\ell(\underline{a}) \leq k \leq \min\{h, \frac{n+\ell(\underline{a})}{2}\}$, the spaces $P_{h,\underline{a}',k}$ are mutually orthogonal.*

PROOF. Each $P_{h,\underline{a}',k}$ is an eigenspace for the self-adjoint operator $D_{h+1,\underline{a}} \times D^*_{h+1,\underline{a}}$. Since the eigenvalue $|Y|(n+\ell(\underline{a})-k-h)(h-k+1)$ is a strictly decreasing function of $k$ for $k \leq \frac{n+\ell(\underline{a})}{2}$, then to different values of $k$ correspond different eigenvalues. This proves the assertion. $\square$

Recall that, if $\underline{a} = (a_0, a_1, \ldots, a_m)$, we set $\underline{a}' = (a_0 - 1, a_1, \ldots, a_m)$ and, inductively, $\underline{a}^{h+1} = \underline{a}^h - (1, 0, \ldots, 0)$.

PROPOSITION 7.11. *Let $F$ be a function in $P_{k,\underline{a}^{h-k},k}$. Then, for $\ell(\underline{a}) \leq k \leq \frac{n+\ell(\underline{a})}{2}$ and $k \leq h \leq n + \ell(\underline{a}) - k$, we have*

$$\|D^*_{h,\underline{a}} D^*_{h-1,\underline{a}'} \cdots D^*_{k+1,\underline{a}^{h-k-1}} F\|^2 = \frac{(n+\ell(\underline{a})-2k)!(h-k)!}{(n+\ell(\underline{a})-k-h)!}|Y|^{h-k}\|F\|^2.$$

*In particular, $D^*_{h,\underline{a}} D^*_{h-1,\underline{a}'} \cdots D^*_{k+1,\underline{a}^{h-k-1}}$ is an isomorphism of $P_{k,\underline{a}^{h-k},k}$ onto $P_{h,\underline{a},k}$.*

PROOF. We prove the assertion by induction on $h$. For $h = k+1$ and $F \in P_{k,\underline{a}',k}$, we have

$$\|D^*_{k+1,\underline{a}} F\|^2 = \langle D^*_{k+1,\underline{a}} F, D^*_{k+1,\underline{a}} F \rangle$$
$$= \langle D_{k+1,\underline{a}} D^*_{k+1,\underline{a}} F, F \rangle$$
$$= |Y|(n+\ell(\underline{a})-2k)\|F\|^2$$

by Proposition 7.8, so the assertion is true. For $h > k+1$, applying Proposition 7.8 to $D_{h,\underline{a}} D^*_{h,\underline{a}}$, we get

$$\|D^*_{h,\underline{a}} D^*_{h-1,\underline{a}'} \cdots D^*_{k+1,\underline{a}^{h-k-1}} F\|^2$$
$$= \langle D_{h,\underline{a}} D^*_{h,\underline{a}} D^*_{h-1,\underline{a}'} \cdots D^*_{k+1,\underline{a}^{h-k-1}} F, D^*_{h-1,\underline{a}'} \cdots D^*_{k+1,\underline{a}^{h-k-1}} F \rangle$$
$$= |Y|(n+\ell(\underline{a})-k-h+1)(h-k)\|D^*_{h-1,\underline{a}'} \cdots D^*_{k+1,\underline{a}^{h-k-1}} F\|^2.$$

Now the proposition follows by induction. $\square$

PROPOSITION 7.12. *Assume $\ell(\underline{a}) \leq h \leq \frac{n+\ell(\underline{a})}{2}$. Then:*

(1) $P_{h,\underline{a}} = \bigoplus_{k=\ell(\underline{a})}^{\min\{h,n+\ell(\underline{a})-h\}} P_{h,\underline{a},k}$;

(2) $D^*_{h+1,\underline{a}} : P_{h,\underline{a}'} \longrightarrow P_{h+1,\underline{a}}$ is an injective map.

PROOF. We prove the assertion by induction on $h$.

Assume that (1) and (2) are true for $\ell(\underline{a}) - 1 \leq h \leq t \leq \frac{n+\ell(\underline{a})-1}{2}$. For $h = \ell(\underline{a}) - 1$ we have $P_{\ell(\underline{a})-1,\underline{a}} = 0$ and so the proposition trivially holds.

Since the operator $D^*_{h,\underline{a}}$ is the adjoint of $D_{h,\underline{a}}$ we have the following decomposition:

$$P_{h,\underline{a}} = Ker(D_{h,\underline{a}}) \oplus D^*_{h,\underline{a}} P_{h-1,\underline{a}'}$$
$$= P_{h,\underline{a},h} \oplus D^*_{h,\underline{a}} P_{h-1,\underline{a}'}.$$

In particular

$$P_{t+1,\underline{a}} = P_{t+1,\underline{a},t+1} \oplus D^*_{t+1,\underline{a}} P_{t,\underline{a}'}.$$

By induction

$$P_{t,\underline{a}'} = \bigoplus_{k=\ell(\underline{a})}^{t} P_{t,\underline{a}',k}$$

and so

$$P_{t+1,\underline{a}} = P_{t+1,\underline{a},t+1} \oplus D^*_{t+1,\underline{a}} \left( \bigoplus_{k=\ell(\underline{a})}^{t} P_{t,\underline{a}',k} \right)$$

$$= \bigoplus_{k=\ell(\underline{a})}^{t+1} P_{t+1,\underline{a},k}.$$

This proves (1), while (2) follows from (1) and Proposition 7.11. $\square$

COROLLARY 7.13. *The dimension of the spaces $P_{h,\underline{a},k}$ that appear in the decomposition of $P_{h,\underline{a}}$ is*

$$\frac{n+\ell(\underline{a})+1-2k}{n-k+1} \binom{n}{k} \binom{k}{\ell(\underline{a})} \binom{\ell(\underline{a})}{a_1,\ldots,a_m} \prod_{j=1}^{m} (\dim(W_j))^{a_j}.$$

PROOF. From the previous proposition it follows that

$$\dim(P_{t+1,\underline{a},t+1}) = \dim(P_{t+1,\underline{a}}) - \dim(P_{t,\underline{a}'}).$$

Now

$$\dim(P_{t+1,\underline{a}}) = \binom{n}{t+1} \binom{t+1}{a_0, a_1, \ldots, a_m} \prod_{j=1}^{m} (\dim(W_j))^{a_j}.$$

In fact, $\binom{n}{t+1}$ represents the number of $(t+1)$-subsets in $X$ and $\binom{t+1}{a_0,a_1,\ldots,a_m} \times \prod_{j=1}^{m}(\dim(W_j))^{a_j}$ represents the number of possible choices in the fundamental function $F = \prod_{r \in A} F^r$ of $a_i$ functions belonging to the eigenspace $W_i$ of $L(Y)$. Thus

$$\dim(P_{t+1,\underline{a},t+1}) = \binom{n}{t+1}\binom{t+1}{a_0,a_1,\ldots,a_m}\prod_{j=1}^{m}(\dim(W_j))^{a_j}$$

$$- \binom{n}{t}\binom{t}{a_0-1,a_1,\ldots,a_m}\prod_{j=1}^{m}(\dim(W_j))^{a_j}$$

$$= \frac{n-t-a_0}{n-t}\binom{n}{t+1}\binom{t+1}{a_0,a_1,\ldots,a_m}\prod_{j=1}^{m}(\dim(W_j))^{a_j}.$$

Since, by Proposition 7.11, $\dim(P_{h,\underline{a},k}) = \dim(P_{k,\underline{a}^{h-k},k})$ one can obtain the result replacing $t$ by $k-1$ and $\underline{a}$ by $\underline{a}^{h-k}$. $\square$

We want to find now the eigenvector of $\frac{\Delta}{norm(\Delta)}$ associated with the eigenvalue 1. Consider in $P_{1,(1,0,\ldots,0)}$ the function

$$f = \sum_{i=1}^{n} f_i,$$

where $f_i$ is the fundamental function of type $(1,0,\ldots,0)$ whose domain is $\{i\}$. Set

$$\langle f \rangle = P_{1,(1,0,\ldots,0),0} =: D^*_{1,(1,0,\ldots,0)} P_{0,(0,\ldots,0),0}.$$

So the element $F_0 = D^*_{h,(h,0,\ldots,0)} \cdots D^*_{3,(3,0,\ldots,0)} D^*_{2,(2,0,\ldots,0)} f$ is the generator of the space $P_{h,(h,0,\ldots,0),0}$, which has dimension 1. Corollary 7.9 implies that $P_{h,(h,0,\ldots,0),0}$ is an eigenspace for $\frac{\Delta}{norm(\Delta)}$ and the corresponding eigenvalue is 1. Moreover, the connectedness of the graph associated with $\Delta$ implies that this is the unique (up to constant) eigenvector of eigenvalue 1. We denote by $P_{1,(1,0,\ldots,0),1}$ the orthogonal subspace to $P_{1,(1,0,\ldots,0),0}$ in $P_{1,(1,0,\ldots,0)}$. It has dimension $n-1$.

Observe that the definition of fundamental functions is strictly linked to the spectral theory of the operator $Q$ and so of the operator $M$ restricted to each domain. In fact, if $F$ is a fundamental function in $P_{h,\underline{a},A}$, with $\underline{a} = (a_0, a_1, \ldots, a_m)$ and $A \in \Omega_h$, then it is an eigenvector for the operator $M$ and the corresponding eigenvalue is $\frac{1}{h}\sum_{j=0}^{m} a_j \lambda_j$. So the set of the eigenvalues of $M$ is given by $\binom{n}{h}$ copies of these values. In particular, the eigenspace $P_{h,\underline{a},k}$ of $\frac{\Delta}{norm(\Delta)}$ is also an eigenspace for $M$ and an eigenvector in this space has eigenvalue $\frac{1}{h}\sum_{j=0}^{m} a_j \lambda_j$. So, by Corollary 7.9 and definition (6) of $P$, we get the following theorem.

THEOREM 7.14. $P_{h,\underline{a},k}$ is an eigenspace for $P$ with eigenvalue

$$p_0 \cdot \frac{1}{h}\sum_{j=0}^{m} a_j \lambda_j + (1-p_0)\frac{(n+\ell(\underline{a})-k-h)(h-k+1)-(n-h)}{h(n-h)}.$$

REMARK 7.15. It is easy to check that the operator $M$ is not ergodic. In fact its associated graph contains $\binom{n}{h}$ connected components and so the multiplicity of the eigenvalue 1 for $M$ is $\binom{n}{h}$.

On the other hand we already observed that the operator $\frac{\Delta}{norm(\Delta)}$ has the eigenvalue 1 with multiplicity 1. To conclude that it is ergodic it suffices to show that $-1$ is not an eigenvalue, that is, the associated graph is not bipartite. In fact consider $\theta \in \Theta_h$ with domain $\{i_1,\ldots,i_h\}$ and $\theta(i_j) = y_j$, for every $j=1,\ldots,h$. By definition of $\Delta$ we can connect $\theta$ with $\varphi$, whose domain is $\{i_1,\ldots,i_{h-1},i_t\}$, $i_h \neq i_t$ and such that $\varphi(i_j) = y_j = \theta(i_j)$ for all $j=1,\ldots,h-1$ and $\varphi(i_t) = y_t$. Moreover $\theta$ can also be connected with $\varrho$ whose domain is $\{i_1,\ldots,i_{h-2},i_h,i_t\}$ and such that $\varrho(i_j) = y_j = \theta(i_j)$ for all $j=1,\ldots,h-2,h$ and $\varrho(i_t) = y_t = \varphi(i_t)$. On the other hand $\varphi$ and $\varrho$ are connected as well and this proves that the graph is not bipartite.

From Proposition 3.4 we can deduce the ergodicity for the operator $P$, since the multiplicity of the eigenvalue 1 is 1 and the eigenvalue $-1$ does not appear in the spectrum of $P$.

REMARK 7.16. The second crested product reduces to the crossed product if $k=n$ and to the nested product if $k=1$.

In fact, if $k=n$, the domain of a function $\theta \in \Theta_n$ cannot be changed and $\theta$ can be identified with the $n$-tuple $(y_1,\ldots,y_n) \in Y^n$ of its images. The operator $P$ becomes

$$P = \frac{1}{n}\sum_{i=1}^{n} I_1 \otimes \cdots \otimes I_{i-1} \otimes Q \otimes I_{i+1} \otimes \cdots \otimes I_n,$$

which is the crossed product on the space $Y^n$.

If $k=1$, then $\Delta$ has the following expression:

$$\Delta = \begin{pmatrix} 0 & 1 & \cdots & & \cdots & 1 \\ 1 & 0 & 1 & & & \vdots \\ \vdots & 1 & \ddots & & & \vdots \\ \vdots & & & \ddots & & 1 \\ 1 & \cdots & & \cdots & 1 & 0 \end{pmatrix}$$

and $norm(\Delta) = n-1$. So we get

$$P = p_0(I_X \otimes Q) + (1-p_0)\left(\frac{\Delta}{norm(\Delta)} \otimes J_Y\right),$$

which is just the nested product of $X$ and $Y$, with $P_X = \frac{\Delta}{norm(\Delta)}$ and $P_Y = Q$.

7.2. *Bi-insect.* In what follows, we take $Y$ as a homogeneous rooted tree of degree $q$ and depth $m-1$ and we give an explicit description of the spectrum of the operator $P = p_0 M + (1-p_0)\frac{\Delta}{norm(\Delta)}$ acting on the space $L(\Theta_2)$. Therefore we are considering functions in $\Theta_2$ such that the image of each element of the domain is an insect. In other words, our aim is to diagonalize the bi-insect Markov chain defined in Section 2.2.2.

Suppose $X$ to be a set of cardinality $n$ and let $m \geq 3$. Recall that we have the decomposition

$$L(Y) = \bigoplus_{j=0}^{m-1} W_j,$$

where $W_0 \cong \mathbb{C}$ and

$$W_j = \left\{ f \in L(L_{m-1}) : f = f(x_1, \ldots, x_j), \sum_{i=0}^{q-1} f(x_1, \ldots, x_{j-1}, i) \equiv 0 \right\},$$

for every $j = 1, \ldots, m-1$. Observe that $\dim(W_j) = q^{j-1}(q-1)$.

The eigenspaces relative to the operator $\Delta/norm(\Delta)$ are the subspaces of the form $P_{2,(a_0,a_1,\ldots,a_{m-1}),k}$, with $k = 0, 1, 2$. The corresponding eigenvalue is

$$\frac{1}{q^{m-1}(n-2)2}[q^{m-1}(n+\ell(\underline{a})-k-2)(2-k+1) - q^{m-1}(n-2)].$$

So, by dependence of $\ell(\underline{a})$, we get the following eigenspaces:

- $P_{2,(a_0,a_1,\ldots,a_{m-1}),2}$ $\begin{cases} a_0 = 0, & \text{with eigenvalue } \lambda = 0, \\ a_0 = 1, & \text{with eigenvalue } \lambda = -\dfrac{1}{2(n-2)}, \\ a_0 = 2, & \text{with eigenvalue } \lambda = -\dfrac{1}{n-2}, \end{cases}$

- $P_{2,(a_0,a_1,\ldots,a_{m-1}),1}$ $\begin{cases} a_0 = 1, & \text{with eigenvalue } \lambda = \dfrac{1}{2}, \\ a_0 = 2, & \text{with eigenvalue } \lambda = \dfrac{n-4}{2(n-2)}, \end{cases}$

- $P_{2,(2,0,\ldots,0),0} \Rightarrow a_0 = 2$   with eigenvalue $\lambda = 1$.

Now we describe the eigenvalues of these eigenspaces with respect to the operator $M$ and to join the results.

If $F$ is a fundamental function of type $(a_0, a_1, \ldots, a_{m-1})$, then it has eigenvalue $\frac{1}{2}\sum_{j=0}^{m-1} a_j \lambda_j$, where $\lambda_j = 1 - \frac{q-1}{q^{m-j}-1}$ is the eigenvalue of the eigenspace $W_j$, of dimension $q^{j-1}(q-1)$, occurring in the spectral decomposition of $L(Y)$. From this we can fill the following table in which we give the eigenspaces, together with the corresponding eigenvalue and dimension.

- $P_{2,(a_0,a_1,\ldots,a_{m-1}),2}$. We have three different cases:

  (1) if $a_0 = 0$, the corresponding eigenspace is

  $$P_{2,(0,\ldots,0,\underbrace{1}_{i\text{th place}},0,\ldots,0,\underbrace{1}_{j\text{th place}},0,\ldots,0),2}$$

  of dimension $n(n-1)(q-1)^2 q^{i-1} q^{j-1}$, with eigenvalue $\frac{p_0}{2}(\lambda_i + \lambda_j)$;

  (2) if $a_0 = 1$, the corresponding eigenspace is

  $$P_{2,(1,\ldots,0,\underbrace{1}_{i\text{th place}},0,\ldots,0),2}$$

  of dimension $n(n-2)(q-1)q^{i-1}$, with eigenvalue $p_0 \frac{1+\lambda_i}{2} + (1-p_0)\frac{-1}{2(n-2)}$;

  (3) if $a_0 = 2$, the corresponding eigenspace is $P_{2,(2,0,\ldots,0),2}$ of dimension $\frac{n(n-3)}{2}$ with eigenvalue $p_0 + (1-p_0)\frac{-1}{n-2}$.

- $P_{2,(a_0,a_1,\ldots,a_{m-1}),1}$. We have two different cases:

  (1) if $a_0 = 1$, the corresponding eigenspace is

  $$P_{2,(1,\ldots,0,\underbrace{1}_{i\text{th place}},0,\ldots,0),1}$$

  of dimension $n(q-1)q^{i-1}$, with eigenvalue $p_0 \frac{1+\lambda_i}{2} + \frac{1-p_0}{2}$;

  (2) if $a_0 = 2$, the corresponding eigenspace is $P_{2,(2,0,\ldots,0),1}$ of dimension $n-1$, with eigenvalue $p_0 + (1-p_0)\frac{n-4}{2(n-2)}$.

- $P_{2,(2,0,\ldots,0),0}$. In this case, the dimension of the eigenspace is 1 with eigenvalue 1.

SECTION DE MATHÉMATIQUES
UNIVERSITÉ DE GENÈVE
2-4, RUE DU LIÈVRE, CASE POSTALE 64
1211 GENÈVE 4
SUISSE
E-MAIL: [Daniele.Dangeli@math.unige.ch](mailto:Daniele.Dangeli@math.unige.ch)
[Alfredo.Donno@math.unige.ch](mailto:Alfredo.Donno@math.unige.ch)